\PassOptionsToPackage{dvipsnames}{xcolor} 
\documentclass[11pt,a4paper]{article}
\usepackage{tikz,amssymb,amsmath,amsthm,multirow,color,xcolor,cancel,bbm}
\usepackage[textheight=25cm,textwidth=16cm]{geometry}
\usepackage{mathrsfs}

\newtheorem{remark}{Remark}
\newtheorem{theorem}{Theorem}
\newtheorem{lemma}{Lemma}
\newtheorem{corollary}{Corollary}
\newtheorem{question}{Question}
\newtheorem{example}{Example}

\def\BF#1{{\boldmath\bf #1}}
\def\cst{{\mathrm{cst}}}
\def\Z{\mathbb{Z}}

\def\h{\tilde h}
\def\Y{\tilde Y}
\def\f{\tilde f}
\def\X{\tilde X}
\def\H{\tilde H}
\def\G{\tilde G}
\def\A{\tilde A}
\def\B{\tilde B}
\def\I{\tilde I}
\def\tpi{\tilde \pi}

\def\din{d^{\mathrm{in}}}
\def\dout{d^{\mathrm{out}}}

\tikzstyle{N}=[outer sep=1,inner sep=2,circle,draw,thick]

\title{Nilpotent dynamics on signed interaction graphs\\ and weak converses of Thomas' rules}

\author{
Adrien Richard
\footnote{Laboratoire I3S, UMR CNRS7271, Universit\'e C\^ote d'Azur, France}
\footnote{CMM,UMI CNRS 2807, Universidad de Chile, Chile}
\footnote{This work is partially supported by the Young Researcher project ANR-18-CE40-0002-01 ``FANs'', and project STIC AmSud CoDANet 19-STIC-03 (Campus France 43478PD).} 
\footnote{E-mail adress: \tt{richard@unice.fr}}
}

\date{January 22, 2018; revised March 28 2019}

\begin{document}

\maketitle

\begin{abstract} 
A finite dynamical system with $n$ components is a function $f:X\to X$ where $X=X_1\times\dots\times X_n$ is a product of $n$ finite intervals of integers. The structure of such a system $f$ is represented by a signed digraph $G$, called interaction graph: there are $n$ vertices, one per component, and the signed arcs describe the positive and negative influences between them. Finite dynamical systems are usual models for gene networks. In this context, it is often assumed that $f$ is {\em degree-bounded}, that is, the size of each $X_i$ is at most the out-degree of $i$ in $G$ plus one. Assuming that $G$ is connected and that $f$ is degree-bounded, we prove the following: if $G$ is not a cycle, then $f^{n+1}$ may be a constant. In that case, $f$ describes a very simple dynamics: a global convergence toward a unique fixed point in $n+1$ iterations. This shows that, in the degree-bounded case, the fact that $f$ describes a complex dynamics {\em cannot} be deduced from its interaction graph. We then widely generalize the above result, obtaining, as immediate consequences, other limits on what can be deduced from the interaction graph only, as the following weak converses of Thomas' rules: if $G$ is connected and has a positive (negative) cycle, then $f$ may have two (no) fixed points. 

\end{abstract}

\section{Introduction}

Let $X=X_1\times \cdots\times X_n$ be a product of $n$ finite intervals of integers, and let $f:X\to X$, 
\[
x=(x_1,\dots,x_n)\mapsto f(x)=(f_1(x),\dots,f_n(x)).
\]
Such a function $f$ is regarded as a {\em finite dynamical system} with $n$ components (FDS for short). If $X=\{0,1\}^n$ then $f$ is usually called a {\em Boolean network} (BN for short). FDSs have many applications. In particular, since the seminal works of MacCulloch and Pitts, Hopfield, Kauffman and Thomas, they are classical models for the dynamics of neural and gene networks \cite{MP43,H82,GM90,K69,K93,T73,TA90}. They are also used in information theory, for the network coding problem \cite{GR11,GRF16}. 

\smallskip
The structure of a FDS $f$ is usually represented by a signed digraph, called {\em interaction graph}. Formally, a {\em signed digraph} is a couple $G=(V,E)$ where $V$ is a set of vertices and $E\subseteq V\times V\times\{+,-\}$ is a set of signed arcs: if $(j,i,s)\in E$ then $G$ has an arc from $j$ to $i$ of sign~$s$. The \BF{interaction graph of $f$} is then defined as the signed digraph  with vertex set $\{1,\dots,n\}$ and with a positive (resp. negative) arc from $j$ to $i$ if there exists $x\in X$ with $x_j<\max(X_j)$ such that 
\[
f_i(x_1,\dots,x_j+1,\dots,x_n)-f_i(x_1,\dots,x_j,\dots,x_n)
\]
is positive (resp. negative). If $G$ is the interaction graph of $f$, we say that $f$ is a \BF{FDS  on $G$}. 

\smallskip
In molecular biology, the first reliable informations available about a biological system often concern its interaction graph \cite{TK01,N15}. A natural question is then the following: {\em what can be said on the dynamics of $f$ according to its interaction graph only?} Among the many dynamical properties that can be studied, fixed points are of special interest, since they correspond to stable states and often have a strong meaning \cite{R86,TA90,GM90}. For instance, in the context of gene networks, they correspond to stable patterns of gene expression at the basis of particular biological processes \cite{TA90}. As such, they are the property which has been studied the most. 

\smallskip 
In particular, many works have been stimulated by two general rules, stated in 1981 by the biologist Ren\'e Thomas \cite{T81}, concerning the influence of positive and negative cycles of the interaction graph (cycles are always directed and without repeated vertices, and the sign of a cycle is the product of the signs of its arcs). Roughly, the first (resp. second) rule says that the presence of a positive (resp. negative) cycle is a necessary condition for multiple stable states (resp. permanent periodic behaviors).  Several variations of these rules have been formally stated and proved for several classes of discrete and continuous dynamical systems (see {\em e.g.} \cite{A08,RRT08,RC07,R08,R09,R18,R19} for the discrete case, and \cite{PM95,S98,G98,CD02,S03,S06,KST07,S13,KS19} for the continuous case). In particular, the first rule has many generalizations and strong versions involving local interaction graphs based on the Jacobian matrix (see \cite{R19,KS19} and the references therein). It has also a nice biological interpretation since multiple stable states can account for cell differentiation processes \cite{TK01a}.  

\smallskip 
In this paper, we are mainly concerned with the basic \BF{discrete versions of Thomas' rules} given by the second and third item of the following theorem. From them we deduce that there is a unique fixed point in the acyclic case, and the first item says something stronger ($f^k$ denotes the $k$th composition of $f$ with itself).    


\begin{theorem}\label{thm:thomas}
Let $f$ a FDS with $n$ components on a signed digraph $G$. 
\begin{itemize}
\item[{\em (a)}] If $G$ has no cycle, then $f^n$ is constant, and thus $f$ has a unique fixed point \cite{R86,R95}. 
\item[{\em (b)}] If $G$ has no negative cycle, then $f$ has at least one fixed point \cite{R10}. 
\item[{\em (c)}] If $G$ has no positive cycle, then $f$ has at most one fixed point \cite{RRT08,RC07}.
\end{itemize}
\end{theorem}

The following complementary theorem shows that, in the the Boolean case, stronger conclusions can be obtained under the additional assumption that the interaction graph is {\em strongly connected} ({\em i.e.} for every distinct vertices $i$ and $j$, there is a path from $i$ to $j$). 

\begin{theorem}\label{thm:aracena} 
Let $f$ be a BN on a strongly connected signed digraph $G$ with at least one arc.
\begin{itemize}
\item[{\em (a)}] If $G$ has no negative cycle, then $f$ has at least two fixed points \cite{A08}. 
\item[{\em (b)}] If $G$ has no positive cycle, then $f$ has no fixed point \cite{A08}.
\end{itemize}
\end{theorem}

The above results give conditions on $G$ that force {\em all} the FDSs or BNs on $G$ to have particular properties. In this paper, we take the converse direction: we study conditions on $G$ that allow {\em at least one} FDS on $G$ to have particular properties. Since it is clearly easier to exhibit a particular FDS when we allow the intervals $X_i$ to be large, it is natural to think to bound the size of these intervals. This leads us to focus on {\em degree-bounded systems}, which are usual in the context of gene networks. The formal definition follows. 

\smallskip
Let $f$ be an FDS with $n$ components, and let $G$ be its interaction graph. We say that $f$ is \BF{degree-bounded} if, for all $i\in \{1,\dots,n\}$, 
\begin{itemize}
\item
$|X_i|=2$ if $\dout_G(i)=0<\din_G(i)$, and 
\item
$|X_i|\leq \dout_G(i)+1$ otherwise. 
\end{itemize}
Here, $\dout_G(i)$ is the out-degree of $i$ in $G$, defined as the  number of positive arcs leaving $i$ plus the number of negative arcs leaving $i$. The in-degree $\din_G(i)$ is defined similarly. 

\smallskip
These conditions on the intervals $X_i$ could seem unnatural at first glance. However, Ren\'e Thomas, who introduced them, showed that they are very natural in the context of gene networks. He actually widely popularized the use of degree-bounded systems for modeling gene networks through the so called {\em generalized logical method}; see \cite{TA90,TK01} for a presentation, and see \cite{AMN14} and the references therein for recent applications. The intuition behind the degree-bounded conditions is the following. In continuous modelings of gene networks, the influence of an interaction from $i$ to $j$ of sign $s$ is often described by a sigmoid function, which can be approximated by a step function with a threshold $\theta^s_{ij}$. From this approximation, the system can be discretized in a natural way: the variable $x_i$, which is a non negative real number representing a concentration, is approximated by a discrete variable $\tilde x_i$ defined as the number of out-going interactions $(i,j,s)$ such that $x_i\geq \theta^s_{ij}$. The size of the domain of this discrete variable $\tilde x_i$ is thus at most the number of interactions leaving $i$ plus one (see \cite{TA90,TK01} for details).   

\subsection{Nilpotent systems}

Our first result concerns {\bf nilpotent} FDSs, that is, FDSs $f$ such that $f^k$ is constant for some positive integer $k$. 

\begin{theorem}\label{thm:nilpotent}
If $G$ is an $n$-vertex connected signed digraph distinct from a signed cycle, then  there exists a degree-bounded system $f$ on $G$ such that $f^{n+1}$ is constant. 
\end{theorem}

This shows that a nilpotent degree-bounded system with a linear convergence can be defined on almost all connected signed digraphs. Since such a nilpotent system describes a very simple dynamics $-$ a global convergence in $n+1$ iterations $-$ this shows that {\em the fact that a degree-bounded system describes a complex dynamics cannot be deduced from its interaction graph.} 

\smallskip
If $G$ is a signed cycle, then any degree-bounded system $f$ on $G$ is a BN, and we deduce from Theorem~\ref{thm:aracena} that $f$ has either zero or two fixed points, thus $f$ is not nilpotent. Therefore, according to the previous theorem, we have the following characterization:

\begin{corollary}\label{cor:nilpotent}
Let $G$ be a connected signed digraph. There exists a nilpotent degree-bounded system on $G$ if and only if $G$ is not a signed cycle. 
\end{corollary}


The existence of nilpotent systems is also the subject of the recent paper \cite{GR16}. It is proved that, for {\em every} $n$-vertex signed digraph $G$,  
\begin{itemize}
\item 
there exists a FDS $f:X\to X$ on $G$ with $X=\{0,1,2,3\}^n$ such that $f^2$ is constant, and
\item
there exists a FDS $f:X\to X$ on $G$ with $X=\{0,1,2\}^n$ such that $f^{\lfloor \log_2 n\rfloor+2}$ is constant. 
\end{itemize}
As shown above, there does not necessarily exist a nilpotent BN on a given signed digraph~$G$, and \cite{GR16} provides several sufficient conditions for the existence of such a BN. This raises the question of a characterization of the signed digraphs on which a nilpotent BN can be defined, in the spirit of the previous corollary, which gives such a characterization in the degree-bounded case. In \cite{GR16} the following question is also raised. 

\begin{question}
Is there a constant $c$ such that, for every $n$-vertex signed digraph $G$, if there exists a nilpotent BN on $G$, then there exists a nilpotent BN $f$ on $G$ that converges in at most $cn$ iterations, that is, such that $f^{cn}$ is constant?
\end{question}

Theorem~\ref{thm:nilpotent} and Corollary \ref{cor:nilpotent} clearly give a positive answer in the degree-bounded case. 

\subsection{Convergence between systems and weak converses of Thomas' rules}

Our second result generalizes the previous one using the following notion of convergence between systems. Let $f:X\to X$ and $h:Y\to Y$ be two FDSs with $n$ components. We say that \BF{$f$ converges toward $h$} in $k$~steps~if 
\begin{itemize}
\item
$f^k(X)\subseteq h(Y)\subseteq Y\subseteq X$, and 
\item
$f(x)=h(x)$ for all $x\in Y$.
\end{itemize} 
In a signed digraph, a {\em source} is a vertex of in-degree zero, a {\em sink} is a vertex of out-degree zero, and an {\em isolated vertex} is a vertex of in- and out-degree zero. 

\begin{theorem}\label{thm:main} 
Let $G$ be a connected signed digraph distinct from a signed cycle. Let $H$ be a subgraph of $G$ obtained by removing arcs only. Suppose that every non-isolated source of $H$ is a source of $G$, and that every non-isolated sink of $H$ is a sink of $G$. Then, for every degree-bounded system $h$ on $H$, there exists a degree-bounded system $f$ on $G$ that converges toward $h$ in $k+1$ steps, where $k$ is the number of isolated vertices in $H$. 
\end{theorem}

Before discussing some consequences of this theorem, let us first observe that some conditions on the subgraph $H$ are needed to obtain the conclusion. Suppose, for instance, that $H$ is obtained from $G$ by removing an arc from $i$ to $j$ with $\din_G(i)=\dout_G(i)=\din_G(j)=1$, and let $h:Y\to Y$ be any degree-bounded system on $H$. Then there is no degree-bounded system $f:X\to X$ on $G$ that converges toward $h$. Indeed, suppose, for a contradiction, that such a $f$ exists. Since $\dout_H(i)< \dout_G(i)=\din_H(i)=1$, $X_i$ and $Y_i$ are of size two, and thus $X_i=Y_i$. Since $f_j$ only depends on $x_i$ and $X_i=Y_i$, there exists $x,y\in Y$ that only differ in $x_i\neq y_i$ such that $f_j(x)\neq f_j(y)$. Since $f$ converges toward $h$, we deduce that $h_j(x)\neq h_j(y)$ and thus $\din_H(j)>0$, a contradiction. 

\smallskip
A first immediate consequence of Theorem~\ref{thm:main} is Theorem~\ref{thm:nilpotent}. Indeed, consider the subgraph $H$ obtained by removing {\em all} the arcs of $G$, and let $h:Y\to Y$ be any degree-bounded system on $H$. Then $H$ consists in $n$ isolated vertices, and thus $Y$ is a singleton, so that $h:\{\xi\}\to\{\xi\}$ for some $\xi\in \Z^n$. Then, according to the conclusion, there exists a degree-bounded system $f$ on $G$ that converges toward $h$ in $n+1$ steps, thus $f^{n+1}=\cst=\xi$, and Theorem~\ref{thm:nilpotent} is recovered. 

\smallskip
As another consequence, we have the following weak converses of Theorems~\ref{thm:thomas}(b) and \ref{thm:thomas}(c). 

\begin{corollary}\label{cor:converse} 
Let $G$ be a connected signed digraph. 
\begin{itemize}
\item[{\em (a)}]
If $G$ has a negative cycle, then there exists a degree-bounded system on $G$ without fixed~point. 
\item[{\em (b)}]
If $G$ has a positive cycle, then there exists a degree-bounded system on $G$ with two fixed~points.
\item[{\em (c)}]
More generally, if $G$ has $k$ vertex-disjoint positive cycles, with $k\geq 1$, then there exists a degree-bounded system on $G$ with $2^k$ fixed points.
\end{itemize}
\end{corollary}

\begin{proof}[{\bf Proof}]
For the first point, suppose that $G$ has a negative cycle $C$, and let $H$ be obtained from $G$ by removing all the arcs of $G$, excepted those of $C$. Let $h$ be any degree-bounded system on $H$. It is easy to see that $h$ has no fixed point. If $G=C$ then $G=H$ and we are done. Otherwise, by Theorem~\ref{thm:main}, there exists a degree-bounded system $f$ on $G$ that converges toward~$h$. Thus $f$ has no fixed point. 

\smallskip
The proof of the third point, which implies the second, is similar. Suppose that $G$ has $k$ vertex-disjoint positive cycles, say $C_1,\dots,C_k$. Let $H$ be obtained from $G$ by removing all the arcs of $G$, excepted those of $C_i$, $1\leq i\leq k$. Let $h$ be any degree-bounded system on $H$. It is easy to see that $h$ has $2^k$ fixed points. If $k=1$ and $G=C_1$, then $G=H$ and we are done. Otherwise, by Theorem~\ref{thm:main}, there exists a degree-bounded system $f$ on $G$ that converges toward~$h$. Thus $f$ has $2^k$ fixed points. 
\end{proof}

Note that from Thomas' rules and the above weak converses, we obtain the following characterizations: {\em if $G$ is connected, then $G$ has a negative (positive) cycle if and only if there exists a degree-bounded system on $G$ no (two) fixed points.}

\smallskip
Note also that the weak converses of Thomas' rules are false in the Boolean case. For instance, if $G$ consists in two cycles with opposite signs that share exactly one vertex, then it is easy to see that {\em every} BN on $G$ has a unique fixed point. More generally, if $G$ is strongly connected, has a unique negative (positive) cycle, and has at least one positive (negative) cycle, then {\em every} BN on $G$ has at least (at most) one fixed point \cite{R18}.

\smallskip
Another gap between the Boolean case and the degree-bounded case, that shows the relevance of Corollary \ref{cor:converse}(c), is the following. Let $\nu^+$ be the maximum number of vertex-disjoint positive cycles in $G$. According to Corollary \ref{cor:converse}(c), the maximum number of fixed points in a degree-bounded system on $G$ is exponential in $\nu^+$ while, in the Boolean case, it can be linear in $\nu^+$. Indeed, as proved in \cite{ARS17}, for every positive integer $k$ there exists a strongly connected signed digraph $G_k$ with $k$ vertex-disjoint positive cycles such that {\em every} BN on $G_k$ has at most $k+1$ fixed points.

\smallskip
The organization of the paper is the following. A slightly stronger version of Theorem~\ref{thm:nilpotent} is proved in Section \ref{sec:nilpotent}, and the proof of Theorem~\ref{thm:main}, which uses this stronger version, is given in Section~\ref{sec:convergence}. 

\section{Nilpotent systems}\label{sec:nilpotent}

\begin{figure}
\[
\begin{array}{c}
\begin{tikzpicture}
\node[N] (1) at (150:1){\scriptsize$1$};
\node[N] (2) at (30:1){\scriptsize$2$};
\node[N] (3) at (270:1){\scriptsize$3$};
\path[->,thick]
(3) edge[Green,bend left=25] (1)
(3) edge[red,bend right=25] (1)
(1) edge[Green] (2)
(2) edge[red] (3)
;
\end{tikzpicture}
\\[3mm]
G
\end{array}
\quad\quad
\begin{array}{c}
\begin{tikzpicture}
\node[N] (1) at (150:1){\scriptsize$1$};
\node[N] (2) at (30:1){\scriptsize$2$};
\node[N] (3) at (270:1){\scriptsize$3$};
\path[->,thick]
(1) edge (2)
(2) edge (3)
(3) edge (1)
;
\end{tikzpicture}
\\[3mm]
|G|
\end{array}
\]
{\caption{\label{fig:pseudo-cycle}
Positive (negative) arcs are green (red). This convention is used throughout the paper. 
}}
\end{figure}
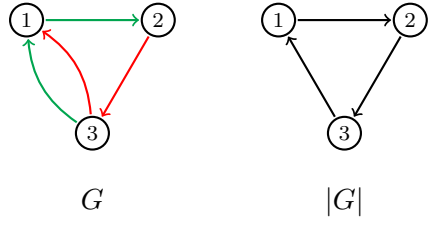

Let $G=(V,E)$ be a signed digraph. The set of {\em positive in-neighbors} of $i$ in $G$, denoted $G^+_i$, is the set of $j\in V$ such that $(j,i,+)\in E$. The set of {\em negative in-neighbors} is defined similarly. The set of in-neighbors of $i$ in $G$ is $G_i=G^+_i\cup G^-_i$. These  notations are rather unconventional, but useful in the following for their compactness. Note that $G^+_i$ and $G^-_i$ are not necessarily disjoint: if $j\in G^+_i\cap G^-_i$ then $G$ has both a positive and a negative arc from $j$ to $i$, and we say that $G$ has {\em parallel arcs} from $j$ to $i$. Note also that $\din_G(i)=|G^+_i|+|G^-_i|$ (thus the in-degree can be greater than the number of in-neighbors and, similarly, the out-degree may be greater than the number of out-neighbors). We denote by $|G|$ be the underlying unsigned digraph of $G$: the vertex set is $V$ and there is an arc from $j$ to $i$ if and only if $j\in G_i$. See Figure~\ref{fig:pseudo-cycle} for an illustration. We say that $G$ is a {\em signed cycle} if $|G|$ is a cycle and $G$ has no parallel arc (thus the signed digraph in Figure~\ref{fig:pseudo-cycle} is not a signed cycle, while its underlying unsigned digraph is a cycle). If $G'=(V',E')$ is another signed digraph, then $G\cup G'=(V\cup V',E\cup E')$. The subgraph of $G$ {\em induced} by $U\subseteq V$, denoted $G[U]$, is the signed digraph with vertex set $U$ and arc set $E\cap U\times U\times \{+,-\}$. We set $G\setminus U=G[V\setminus U]$. We say that $U$ is a {\em strong component} of $G$ is $G[U]$ is strongly connected and $U$ is maximal for this property. A strong component $U$ is {\em initial} if $G$ has no arc from $V\setminus U$ to $U$. A strong component $U$ is {\em trivial} if $G[U]$ is the {\em trivial graph}, {\em i.e.} has one vertex and no arc. 

\smallskip
We say that $G$ is {\bf basic} is all the initial strong components of $G$ are trivial. We set 
\[
\beta(G)=
\begin{cases}
0&\text{if $G$ is basic},\\
1&\text{otherwise.}
\end{cases}
\]
If $i,j\in V$ then $d_G(j,i)$ denotes the minimum number of arcs in a path from $j$ to $i$. If $j=i$ then $d_G(j,i)=0$ and if there is no path from $j$ to $i$ then $d_G(j,i)=\infty$. If $I\subseteq V$ then $d(I,i)$ is the minimum of $d(j,i)$ for $j\in I$. Let $\mathcal{I}$ be the set of initial strong components of $G$. We set 
\[
\lambda(G)=\max_{i\in V}\min_{I\in\mathcal{I}} d_G(I,i)+|I|.
\]

\begin{remark}\label{rem:lambda}
{\em Suppose that $G$ has $n$ vertices. Clearly, for every vertex $i$, there exists $I^*\in\mathcal{I}$ such that $G$ has a path from $I^*$ to $i$. If $P$ is a shortest path from $I^*$ to $i$, then, excepted the initial vertex of $P$, no vertex belong to an initial strong component of $G$. We deduce that the number of arcs in $P$ is at most $n-\sum_{I\in\mathcal{I}}|I|$. As a consequence, 
\[
\lambda(G)\leq n-\sum_{I\in\mathcal{I}}|I|+\max_{I\in\mathcal{I}} |I|\leq n.
\]
In particular, if $G$ has at least two initial strong components then $\lambda(G)< n$, and if $G$ is basic and has $k$ sources, then $\lambda(G)\leq n-k+1$. Note also that $\lambda(G)=n$ if $G$ is strongly connected.}
\end{remark}
 
We are now in position to state the following quantitative version of Theorem~\ref{thm:nilpotent}. We write $f^k=\cst=\xi$ to mean that $f^k$ is a constant function, always equal to $\xi$.  

\begin{theorem}\label{thm:nilpotent2}
Let $G$ be a connected signed digraph distinct from a signed cycle. There exists a degree-bounded system $f:X\to X$ on $G$ such that 
\[
f^{\lambda(G)+\beta(G)}=\cst=\xi, 
\]
where $\xi_i=\min(X_i)$ for all sources $i$ of $G$. 
\end{theorem}

\begin{proof}[{\bf Proof}]
Suppose that $G$ has vertex set $V=\{1,\dots,n\}$. If $G$ is trivial the theorem is obvious, so assume that $G$ is not trivial. The proof is devised in two cases.    

\smallskip
{\bf Case 1: {\boldmath$|G|$\unboldmath} is not a cycle.} Let $\mathcal{I}=\{I_1,\dots,I_k\}$ be the set of initial strong components of $G$. Since $G$ is connected and since $|G|$ is not a cycle, each $I_\ell$ contains at least one vertex $i_\ell$ such that either $i_\ell$ is a source or, for all $j\in G_{i_\ell}$, the out-degree of $j$ in $|G|$ is at least two.  Let $H$ be the subgraph obtained from $G$ by removing all arcs with $i_\ell$ as terminal vertex, $\ell=1,\dots,k$. The initial strong components of $H$ are $\{i_1\},\dots,\{i_k\}$ thus $H$ is basic. Furthermore, for every $\ell=1,\dots,k$ and $j\in G_{i_\ell}$, the out degree of $j$ in $H$ is at least one. Hence, 
\begin{equation}\tag{$*$}\label{eq:sink}
\text{\em every sink of $H$ is a sink of $G$}. 
\end{equation}
Furthermore,
\[
\lambda(H)\leq\lambda(G). 
\]
Indeed, let $i\in V$ and let any $\ell$ such that $G$ has a path $P$ from $I_\ell$ to $i$. Suppose that $P$ is of length  $d_G(I_\ell,i)$. Let $j$ be the first vertex of $P$, and let $P'$ be a path from $i_\ell$ to $j$ of length $d_G(i_\ell,j)$. Since $P'\cup P$ is a path of $H$, $d_H(i_\ell,i)\leq d_G(i_\ell,j)+d_G(I_\ell,i)$, and since $P'$ is contained in $G[I_\ell]$ we have $d_G(i_\ell,j)<|I_\ell|$. Hence $d_H(i_\ell,i)+1\leq d_G(I_\ell,i)+|I_\ell|$. It follows that $\lambda(H)\leq \lambda(G)$. 

\smallskip
Let $I=\{i_1,\dots,i_k\}$ and let $X=X_1\times\dots\times X_n$ be defined as follows: for all $i\in V$,
\[
X_i=
\left\{
\begin{array}{lll}
\{0,1,2,3\}&\text{if $i\in G^+_j\cap G^-_j$ for some $j\in I$}&(1)\\
\{0,1,2\}&\text{if $i\in G^+_j\setminus G^-_j$ for some $j\in I$}&(2)\\
\{0,1,2\}&\text{if $i\in G^-_j\setminus G^+_j$ for some $j\in I$}&(3)\\
\{0,1,2\}&\text{if $i\in G^+_j\cap G^-_j$ for some $j\in V\setminus I$}&(4)\\
\{0,1\}&\text{otherwise.}&(5)
\end{array}
\right.
\]
Let $f:X\to X$ be any FDS on $G$. Since $G$ has no isolated vertex, to prove that $f$ is degree-bounded, it is sufficient to prove that $|X_i|\leq\max(2,\dout_G(i)+1)$ for every $i$. This is obvious in cases $(4)$ and $(5)$. In case (1) we have $\dout_G(i)\geq 2$, and if $\dout_G(i)=2$, then $i$ is a sink of $H$ but not a sink of $G$, a contradiction with (\ref{eq:sink}). Hence $\dout_G(i)\geq 3$ and thus $|X_i|=4\leq \dout_G(i)+1$. Similarly, in case (2) and (3) we have $\dout_G(i)\geq 1$, and if $\dout_G(i)=1$, then $i$ is a sink of $H$ but not a sink of $G$, a contradiction with (\ref{eq:sink}). Hence $\dout_G(i)\geq 2$ and thus $|X_i|=3\leq \dout_G(i)+1$. Thus every FDS $f:X\to X$ on $G$ is degree-bounded. 

\smallskip
For $1\leq p\leq \lambda(H)$, let $J_p$ be the set of $i\in V$ such that $d_H(I,i)+1=p$. In this way, $J_1=I$, and for all $i\in J_p$ with $p\geq 2$, $G$ has at least on arc from $J_{p-1}$ to $i$. Furthermore, since $\lambda(H)=\max_{i\in V}d_H(I,i)+1$, every vertex $i\in V$ is contained in some $J_p$. 

\smallskip
Let $\xi\in X$ be inductively defined as follows: 
\begin{itemize}
\item
For all $i\in J_1$:
\[
\xi_i=
\left\{
\begin{array}{ll}
1&\text{if $G^-_i\setminus G^+_i\neq\emptyset$}\\[1mm]
0&\text{otherwise.}
\end{array}
\right.
\]
\item
For all $i\in J_p$ with $p\geq 2$:
\[
\xi_i=
\left\{
\begin{array}{ll}
1&\text{if $\xi_j=1$ for all $j\in G^+_i\cap J_{p-1}$}\\
&\quad\text{and $\xi_j=0$ for all $j\in (G^-_i\setminus G^+_i)\cap J_{p-1}$}\\[1mm]
0&\text{otherwise.}
\end{array}
\right.
\]
\end{itemize}

\smallskip
Finally, let $f:X\to X$ be defined as follows:
\begin{itemize}
\item
For all $i\in J_1$:
\[
f_i(x)=
\left\{
\begin{array}{ll}
1&\text{if $x_j\geq 2$ for some $j\in G^+_i\setminus G^-_i$}\\
 &\quad\text{or $x_j=2$ for some $j\in G^+_i\cap G^-_i$}\\
 &\quad\quad\text{or $x_j<2$ for some $j\in G^-_i\setminus G^+_i$}\\[1mm]
0&\text{otherwise.}
\end{array}
\right.
\]
\item
For all $i\in V\setminus J_1$ with $\xi_i=1$:
\[
f_i(x)=
\left\{
\begin{array}{ll}
1&\text{if $x_j\geq 1$ for some $j\in G^+_i\setminus G^-_i$}\\
&\quad\text{or $x_j=1$ for some $j\in G^+_i\cap G^-_i$}\\
&\quad\quad\text{or $x_j<1$ for some $j\in G^-_i\setminus G^+_i$}\\[1mm]
0&\text{otherwise.}
\end{array}
\right.
\]
\item
For all $i\in V\setminus J_1$ with $\xi_i=0$:
\[
f_i(x)=
\left\{
\begin{array}{ll}
1&\text{if $x_j\geq 1$ for all $j\in G^+_i\setminus G^-_i$}\\
&\quad\text{and $x_j=1$ for all $j\in G^+_i\cap G^-_i$}\\
&\quad\quad\text{and $x_j<1$ for all $j\in G^-_i\setminus G^+_i$}\\[1mm]
0&\text{otherwise.}
\end{array}
\right.
\]
\end{itemize}

\smallskip
It is straightforward to prove that $G$ is the interaction graph of $f$. Let $\beta=\beta(G)$. Since $\lambda(H)\leq\lambda(G)$, and according to the definition of $\xi$, to conclude the proof it is sufficient to prove that $f^{\lambda(H)+\beta}=\cst=\xi$. For that we prove the following:
\[\tag{$**$}\label{eq:induction}
\forall 1\leq p\leq\lambda(H),~\forall i\in J_p,~\forall \ell\geq p,\qquad 
f^{\ell+\beta}_i=\cst=\xi_i
\]
We proceed by induction on $p$. 
\begin{itemize}
\item 
Let $i\in J_1$. If $\beta=0$ then $i$ is a source of $G$, thus $G^-_i\setminus G^+_i=\emptyset$ and we deduce that $f_i=\cst=0=\xi_i$, thus $f^\ell_i=\cst=\xi_i$ for all $\ell\geq 1$. Suppose now that $\beta=1$. Let $x\in X$ and $\ell\geq 1$. Since $f(X)\subseteq\{0,1\}^n$ we have $f^\ell(x)\in\{0,1\}^n$. If $\xi_i=1$ there exists $j\in G^-_i\setminus G^+_i$ and since $f^\ell_j(x)<2$ we have $f_i(f^\ell(x))=1=\xi_i$. If $\xi_i=0$ then $G^-_i\setminus G^+_i=\emptyset$, and since $f^\ell_j(x)<2$ for all $j\in G^+_i$, we deduce that $f_i(f^\ell(x))=0=\xi_i$. Thus $f_i^{\ell+1}=\cst=\xi_i$ in all cases. 
\item
Suppose that $i\in J_p$ with $p\geq 2$. Let $x\in X$ and $\ell\geq p$. If $\xi_i=1$ then by the definition of $\xi$ and the induction hypothesis, we have $f^{\ell-1+\beta}_j(x)=\xi_j=1$ for all $j\in G^+_i\cap J_{p-1}$ and $f^{\ell-1+\beta}_j(x)=\xi_j=0$ for all $j\in (G^-_i\setminus G^+_i)\cap J_{p-1}$. Since $(G^+_i\cup G^-_i)\cap J_{p-1}$ is not empty, we deduce from the definition of $f_i$ that $f_i(f^{\ell-1+\beta}(x))=1=\xi_i$. If $\xi_i=0$ then by the definition of $\xi$ and the induction hypothesis, we have $f^{\ell'}_j(x)=\xi_j=0$ for some $j\in G^+_i\cap J_{p-1}$ or $f^{\ell'}_j(x)=\xi_j=1$ for some $j\in (G^+_i\setminus G^-_i)\cap J_{p-1}$, and we deduce from the definition of $f_i$ that $f_i(f^{\ell-1+\beta}(x))=0=\xi_i$. Thus $f_i^{\ell+\beta}=\cst=\xi_i$ in all cases.
\end{itemize}     
This proves (\ref{eq:induction}) and completes the proof of the first case. See Example~\ref{ex1} for an illustration. 

\smallskip
{\bf\boldmath Case 2: $|G|$ is a cycle.\unboldmath} Then $\lambda(G)=n$ and $\beta(G)=1$, and since $G$ has no source, we have to prove that there exists a degree-bounded $f$ system on $G$ such that $f^{n+1}=\cst$. Without loss of generality, suppose that $|G|$ has an arc from $i$ to $i+1$ for all $1\leq i<n$ and an arc from $n$ to $1$. For convenience, identify $n+1$ with $1$ and $0$ with $n$. Let $I$ be the set of vertices $i$ such that $G$ has parallel arcs from $i$ to $i+1$. Since $G$ is not a signed cycle and since $|G|$ is a cycle (as in Figure~\ref{fig:pseudo-cycle}), $I$ is not empty, and without loss of generality we assume that $n\in I$. 

\smallskip
Let $X=X_1\times\dots\times X_n$ be defined by 
\[
X_i=
\begin{cases}
\{0,1,2\}&\text{if }i\in I\\
\{0,1\}&\text{otherwise.}
\end{cases}
\]
Let $f:X\to X$ be defined by
\[
f_i(x)=
\left\{
\begin{array}{ll}
\max(X_i)&\text{if $x_{i-1}=1$ and $i-1\in G^+_i$}\\
&\quad\text{or $x_{i-1}=0$ and $i-1\in G^-_i\setminus G^+_i$}\\[1mm]
0&\text{otherwise.}
\end{array}
\right.
\]

\smallskip
It is straightforward to show that $f$ is a degree-bounded system on $G$. Let us prove that $f^{n+1}=\cst$. Since $f_i$ only depends on $x_{i-1}$, we can abusively write $f_i(x)=f_i(x_{i-1})$; this allows us to define $\xi\in X$ by induction in the following way: $\xi_1=0$ and $\xi_i=f_i(\xi_{i-1})$ for all $2\leq i\leq n$. Let us prove, by induction on $i$, that $f^{\ell+1}_i(x)=\xi_i$ for $x\in X$ and $\ell\geq i$. Let $x\in X$. Since $n\in I$, we have $n\in G^+_1\cap G^-_1$ and $f^\ell_n(x)\in\{0,2\}$ for all $\ell\geq 1$, thus $f_1(f^\ell(x))=0=\xi_1$ and the base case is proved. Let $1<i\leq n$ and $\ell\geq i$. By induction hypothesis $f^\ell_{i-1}(x)=\xi_{i-1}$ thus $f_i(f^\ell(x))=f_i(\xi_{i-1})=\xi_i$ and the induction step is proved. Thus $f^{n+1}=\cst=\xi$.
\end{proof} 

\begin{remark}
{\em The bound $\lambda(G)+\beta(G)$ can be reached. For instance, if $G$ satisfies one of the following conditions, then it is straightforward to prove that there is no degree-bounded system $f$ on $G$ such that $f^{\lambda(G)+\beta(G)-1}$ is constant:
\begin{itemize}
\item
$G$ is acyclic,
\item
$G$ is obtained from an acyclic signed digraph by adding a loop on each source,
\item
$G$ consists of two cycles that intersect in a path.
\end{itemize}
It could be interesting to characterize the signed digraphs $G$ for which the bound is reached.}
\end{remark}

\begin{remark}\label{rm:component}
{\em Let $G$ be a signed digraph, and let $G_1,\dots,G_k$ be its connected components. It is easy to see that $\lambda(G)=\max_\ell\lambda(G_\ell)$, and that $\beta(G)=\max_\ell \beta(G_\ell)$. Hence, if $G_\ell$ is not a signed cycle for all $\ell$, by applying the previous theorem on each connected component, we deduce that there exists a degree-bounded system $f:X\to X$ on $G$ such that $f^{\lambda(G)+\beta(G)}=\cst=\xi$.}
\end{remark}

\begin{example}\label{ex1}
{\em Let $G$ be as follows:
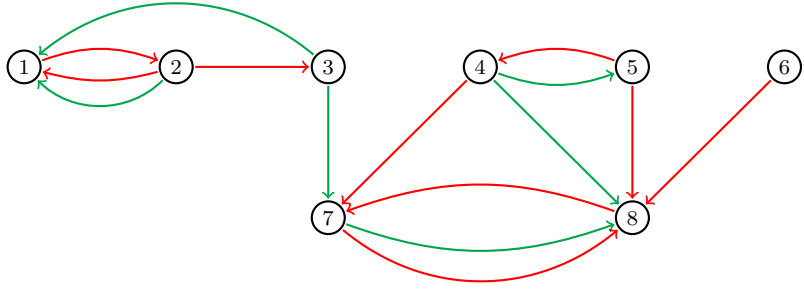
\[
\begin{tikzpicture}
\node[outer sep=1,inner sep=2,circle,draw,thick] (1) at (0,0){\scriptsize$1$};
\node[outer sep=1,inner sep=2,circle,draw,thick] (2) at (2,0){\scriptsize$2$};
\node[outer sep=1,inner sep=2,circle,draw,thick] (3) at (4,0){\scriptsize$3$};
\node[outer sep=1,inner sep=2,circle,draw,thick] (4) at (6,0){\scriptsize$4$};
\node[outer sep=1,inner sep=2,circle,draw,thick] (5) at (8,0){\scriptsize$5$};
\node[outer sep=1,inner sep=2,circle,draw,thick] (6) at (10,0){\scriptsize$6$};
\node[outer sep=1,inner sep=2,circle,draw,thick] (7) at (4,-2){\scriptsize$7$};
\node[outer sep=1,inner sep=2,circle,draw,thick] (8) at (8,-2){\scriptsize$8$};
\path[->,thick]
(1) edge[red,bend left=20] (2)
(2) edge[red,bend left=15] (1)
(2) edge[Green,bend left=45] (1)
(3) edge[Green,bend right=40] (1)
(2) edge[red] (3)
(4) edge[Green,bend right=20] (5)
(5) edge[red,bend right=20] (4)
(3) edge[Green] (7)
(4) edge[red] (7)
(4) edge[Green] (8)
(5) edge[red] (8)
(6) edge[red] (8)
(7) edge[Green, bend right=20] (8)
(7) edge[red, bend right=40] (8)
(8) edge[red, bend right=20] (7)
;
\end{tikzpicture}
\]
The initial strong components are $I_1=\{1,2,3\}$, $I_2=\{4,5\}$ and $I_3=\{6\}$. We have $\lambda(G)=3$, and since $I_1$ and $I_2$ are not trivial, we have $\beta(G)=1$. Then, for the construction described in the proof, we must take $i_1=1$ and $i_3=6$, while $i_2$ should be $4$ or $5$. Let us take $i_2=4$. Then $I=\{1,4,6\}$ and $X$ is as follows:
\[
\begin{array}{l}
X_1=\{0,1\}\\
X_2=\{0,1,2,3\}\\
\end{array}
\qquad
\begin{array}{l}
X_3=\{0,1,2\}\\
X_4=\{0,1\}
\end{array}
\qquad
\begin{array}{l}
X_5=\{0,1,2\}\\
X_6=\{0,1\}\\
\end{array}
\qquad
\begin{array}{l}
X_7=\{0,1,2\}\\
X_8=\{0,1\}\\
\end{array}
\]
Besides, $H$ is obtained by removing the parallel arcs from $2$ to $1$, the arc from $3$ to $1$, and the arc from $5$ to $4$. Thus $\lambda(H)=3$ with $J_1=I=\{1,4,6\}$, $J_2=\{2,5,7,8\}$, and $J_3=\{3\}$. We deduce that 
\[
\begin{array}{l}
\xi_1=0\\
\xi_4=1\\
\xi_6=0\\
\end{array}
\qquad
\begin{array}{l}
\xi_2=1\\
\xi_5=1\\
\xi_7=0\\
\xi_8=1
\end{array}
\qquad
\begin{array}{l}
\xi_3=0.
\end{array}
\]
and thus $f$ is defined by
\[
\begin{array}{l}
f_1(x)=
\begin{cases}
1\text{ if $x_2=2$}\\
\text{~~~ or $x_3\geq 2$,}\\
0\text{ otherwise}
\end{cases}
\\[8mm]
f_4(x)=
\begin{cases}
1\text{ if $x_5<2$,}\\
0\text{ otherwise}
\end{cases}
\\[6mm]
f_6(x)=0
\end{array}
~
\begin{array}{l}
f_2(x)=
\begin{cases}
1\text{ if $x_1<1$}\\
0\text{ otherwise}
\end{cases}
\\[5mm]
f_5(x)=
\begin{cases}
1\text{ if $x_4\geq 1$,}\\
0\text{ otherwise}
\end{cases}
\\[5mm]
f_7(x)=
\begin{cases}
1\text{ if $x_3\geq 1$}\\
\text{~~~ and $x_4<1$,}\\
0\text{ otherwise}
\end{cases}
\\[8mm]
f_8(x)=
\begin{cases}
1\text{ if $x_4\geq 1$}\\
\text{~~~ or $x_5<1$}\\
\text{~~~~~~ or $x_6<1$}\\
\text{~~~~~~~~~ or $x_7=1$,}\\
0\text{ otherwise}
\end{cases}
\end{array}
\!\!
\begin{array}{l}
f_3(x)=
\begin{cases}
1\text{ if $x_2<1$,}\\
0\text{ otherwise.}
\end{cases}
\end{array}
\]
It is easy to check that the interaction graph of $f$ is $G$. Furthermore, if $x\in\{0,1\}^n$, then it is easy to see that, for all $\ell\geq 1$,
\[
\begin{array}{l}
f^\ell_1(x)=0=\xi_1\\
f^\ell_4(x)=1=\xi_2\\
f^\ell_6(x)=0=\xi_6\\
\end{array}
\qquad
\begin{array}{l}
f^{\ell+1}_2(x)=1=\xi_2\\
f^{\ell+1}_3(x)=1=\xi_3\\
f^{\ell+1}_7(x)=0=\xi_7\\
f^{\ell+1}_8(x)=1=\xi_8\\
\end{array}
\qquad
\begin{array}{l}
f^{\ell+2}_3(x)=0=\xi_3.
\end{array}
\]
Since $f(x)\in\{0,1\}^n$ for all $x\in X$, we deduce that $f^4(x)=\xi$ for all $x\in X$, as desired.}
\end{example}

\begin{example}
{\em Suppose that $G$ is the $3$-vertex signed digraph as in Figure~\ref{fig:pseudo-cycle}:
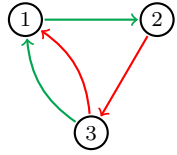
\[
\begin{tikzpicture}
\node[outer sep=1,inner sep=2,circle,draw,thick] (1) at (150:1){\scriptsize$1$};
\node[outer sep=1,inner sep=2,circle,draw,thick] (2) at (30:1){\scriptsize$2$};
\node[outer sep=1,inner sep=2,circle,draw,thick] (3) at (270:1){\scriptsize$3$};
\path[->,thick]
(3) edge[Green,bend left=25] (1)
(3) edge[red,bend right=25] (1)
(1) edge[Green] (2)
(2) edge[red] (3)
;
\end{tikzpicture}
\]
Then $|G|$ is a cycle, and we fall in the second case of the proof. We then have $I=\{3\}$, thus $X_1=X_2=\{0,1\}$ and $X_3=\{0,1,2\}$. Therefore, $f$ is defined by 
\[
f_1(x)=
\begin{cases}
1\text{ if $x_3=1$}\\
0\text{ otherwise}
\end{cases}
\quad
f_2(x)=
\begin{cases}
1\text{ if $x_1=1$}\\
0\text{ otherwise}
\end{cases}
f_3(x)=
\begin{cases}
2\text{ if $x_2=0$}\\
0\text{ otherwise.}
\end{cases}
\]
It is easy to check that the interaction graph of $f$ is $G$. Furthermore, $f_3(x)\neq 1$ for all $x\in X$. Thus, for all $\ell\geq 2$, we have $f^\ell_1(x)=0$ which implies $f^{\ell+1}_2(x)=0$, which implies $f^{\ell+2}_3(x)=2$. Thus $f^4(x)=(0,0,2)$ for all $x\in X$, as desired.} 
\end{example}

\begin{example}
{\em Suppose that $G$ is the $1$-vertex signed digraph as follows:
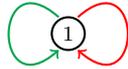
\[
\begin{tikzpicture}
\node[outer sep=1,inner sep=2,circle,draw,thick] (1) at (0,0){\scriptsize$1$};
\draw[Green,->,thick] (1.120) .. controls (-1,1) and (-1,-1) .. (1.-120);
\draw[red,->,thick] (1.60) .. controls (1,1) and (1,-1) .. (1.-60);
\end{tikzpicture}
\]
Then $|G|$ is a cycle, and we fall in the second case of the proof. We then have $I=\{1\}$, thus $X=X_1=\{0,1,2\}$, and $f=f_1$ is defined by $f(0)=0$, $f(1)=2$ and $f(2)=0$. Thus $f^2(x)=0$ for all $x\in X$, as desired.}
\end{example}

\section{Convergence between systems}\label{sec:convergence}

In this section we prove the following strengthening of Theorem~\ref{thm:main}, which is more suited for a proof by induction. A {\em spanning} subgraph is a subgraph obtained by removing arcs only. 

\begin{theorem}\label{thm:main2} 
Let $G$ be signed digraph and let $H$ be a spanning subgraph of $G$. Let $I$ be the set of vertices that are isolated in $H$ but not in $G$. Suppose that every source of $H\setminus I$ is a source of $G$, and that every sink of $H\setminus I$ is a sink of $G$. Suppose also that no connected component of $G$ is a signed cycle with only vertices in $I$. Then, for every degree-bounded system $h$ on $H$ there exists a degree-bounded system $f$ on $G$ that converges toward $h$ in at most $|I|+1$ steps.
\end{theorem}

The proof of Theorem~\ref{thm:main2} uses Theorem~\ref{thm:nilpotent2} and the following technical lemma.

\begin{lemma}\label{lem:extension} 
Let $G$ be a signed digraph and let $H$ be a spanning subgraph of $G$ such that
\begin{itemize}
\item
$G$ has no arc from a sink of $H$ to a source of~$H$, and
\item
$G$ has no arc from a vertex of $H$ to an isolated vertex of $H$. 
\end{itemize}
Let $A$ be the set of sources of $H$ that are not sources of $G$, and let $B$ be the set sinks of $H$ that are not sinks of $G$. Let $h:Y\to Y$ be a degree-bounded system on $H$ and $\xi\in Y$. There exists a degree-bounded system $f:X\to X$ on $G$ such that
\begin{itemize}
\item 
$f(X)\subseteq Y$,
\item
$f_i(X)\subseteq h_i(Y)$ for all $i\not\in A$,
\item
$f(x)=h(x)$ for all $x\in Y$ such that $x_i=\xi_i$ for all $i\in B$, and
\item
$f_i(x)=h_i(x)$ for all $x\in Y$ and $i$ such that $G_i\cap B=\emptyset$.
\end{itemize}
\end{lemma}

\begin{proof}[{\bf Poof of Theorem~\ref{thm:main2}, assuming Lemma~\ref{lem:extension}}] 
Without loss of generality, suppose that $V=\{1,\dots,n\}$ is the vertex set of $G$ and  that $I=\{1,\dots,m\}$ for some $1\leq m\leq n$. Let $h:Y\to Y$ be a degree-bounded system on $H$. Since $h$ is degree-bounded, for all $i\in I$ there exists $\xi_i\in\Z$ such that $Y_i=\{\xi_i\}$.   

\smallskip
Suppose first that $I$ is empty. Then $G$ and $H$ have the same isolated vertices and thus the same sources and sinks. Hence, the conditions of Lemma \ref{lem:extension} are satisfied, and the sets $A$ and $B$ of the statement are both empty. Thus there exists a degree-bounded system $f$ on $G$ that converges toward $h$ in one step, and thus $f$ has the desired properties. Hence, in the following, we assume that $I$ is not empty.

\smallskip
We say that $I$ has the \BF{property $P$} if every connected component $F$ of $G[I]$ satisfies at least one of the following conditions:
\begin{itemize}
\item $F$ is not strongly connected.
\item $F$ has a {\em leaving arc}, that is, $G$ has an arc from a vertex of $F$ to $V\setminus I$.
\item $F$ has no {\em entering arc}, that is, $G$ has no arc from $V\setminus I$ to a vertex of $F$.
\end{itemize}

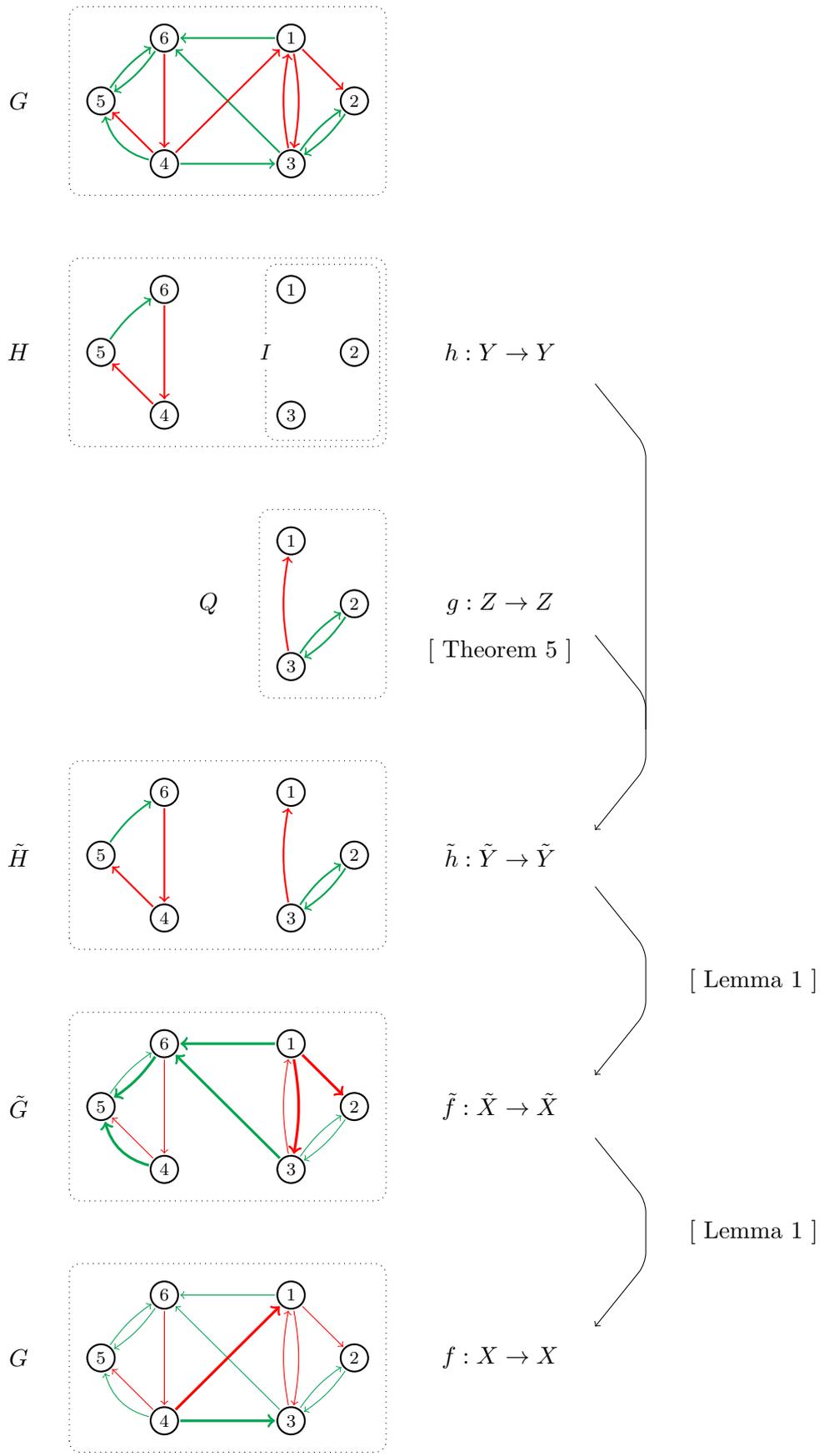
\begin{figure}
\[
\begin{tikzpicture}
\draw[dotted,rounded corners=5] (-0.5,-0.5) rectangle (4.5,2.5);
\node at (-1.3,1) {$G$};
\node[N] (1) at (1,2){\scriptsize$6$};
\node[N] (2) at (3,2){\scriptsize$1$};
\node[N] (3) at (4,1){\scriptsize$2$};
\node[N] (4) at (3,0){\scriptsize$3$};
\node[N] (5) at (1,0){\scriptsize$4$};
\node[N] (6) at (0,1){\scriptsize$5$};
\path[->,thick]
(1) edge[red] (5)
(5) edge[red] (6)
(5) edge[Green,bend left=30] (6)
(6) edge[Green,bend left=10] (1)
(1) edge[Green, bend left=10] (6)
(4) edge[red,bend left=10] (2)
(2) edge[red] (3)
(2) edge[red,bend left=10] (4)
(3) edge[Green,bend left=10] (4)
(4) edge[Green,bend left=10] (3)
(2) edge[Green] (1)
(4) edge[Green] (1)
(5) edge[red](2)
(5) edge[Green](4)
;
\begin{scope}[shift={(0,-4)}]
\draw[dotted,rounded corners=5] (-0.5,-0.5) rectangle (4.5,2.5);
\node at (-1.3,1) {$H$};
\draw[dotted,rounded corners=5] ({3-0.4},{0-0.4}) rectangle ({4+0.4},{2+0.4});
\node[preaction={fill,white}] at ({3-0.4},1) {\small $I$};
\node at (6.3,1) {$h:Y\to Y$};
\node[N] (1) at (1,2){\scriptsize$6$};
\node[N] (2) at (3,2){\scriptsize$1$};
\node[N] (3) at (4,1){\scriptsize$2$};
\node[N] (4) at (3,0){\scriptsize$3$};
\node[N] (5) at (1,0){\scriptsize$4$};
\node[N] (6) at (0,1){\scriptsize$5$};
\path[->,thick]
(1) edge[red] (5)
(5) edge[red] (6)
(6) edge[Green,bend left=10] (1)
;
\draw[rounded corners=5,->] 
(7.8,{1-0.5}) -- (8.6,{1-1.5}) -- (8.6,{-7.1+1.5}) -- (7.8,{-7.1+0.5})
;
\end{scope}
\begin{scope}[shift={(0,-8)}]
\draw[dotted,rounded corners=5] (2.5,-0.5) rectangle (4.5,2.5);
\node at (1.7,1) {$Q$};
\node at (6.3,1) {$g:Z\to Z$};
\node at (6.3,0.25) {[~Theorem~\ref{thm:nilpotent2}~]};
\node[N] (2) at (3,2){\scriptsize$1$};
\node[N] (3) at (4,1){\scriptsize$2$};
\node[N] (4) at (3,0){\scriptsize$3$};
\path[thick,->]
(4) edge[red,bend left=10] (2)
(3) edge[Green,bend left=11] (4)
(4) edge[Green,bend left=11] (3)
;
\draw[rounded corners=5,-] 
(7.8,{1-0.5}) -- (8.6,{1-1.5}) -- (8.6,-1)
;\end{scope}
\begin{scope}[shift={(0,-12)}]
\draw[dotted,rounded corners=5] (-0.5,-0.5) rectangle (4.5,2.5);
\node at (-1.3,1) {$\H$};
\node at (6.3,1) {$\h:\Y\to \Y$};
\node[N] (1) at (1,2){\scriptsize$6$};
\node[N] (2) at (3,2){\scriptsize$1$};
\node[N] (3) at (4,1){\scriptsize$2$};
\node[N] (4) at (3,0){\scriptsize$3$};
\node[N] (5) at (1,0){\scriptsize$4$};
\node[N] (6) at (0,1){\scriptsize$5$};
\path[thick,->]
(1) edge[red] (5)
(5) edge[red] (6)
(6) edge[Green,bend left=10] (1)
(4) edge[red,bend left=10] (2)
(3) edge[Green,bend left=11] (4)
(4) edge[Green,bend left=11] (3)
;
\draw[rounded corners=5,->] 
(7.8,{1-0.5}) -- (8.6,{1-1.5}) -- (8.6,{-3+1.5}) -- (7.8,{-3+0.5})
;
\node at (10.3,-1) {[~Lemma~\ref{lem:extension}~]};
\end{scope}
\begin{scope}[shift={(0,-16)}]
\draw[dotted,rounded corners=5] (-0.5,-0.5) rectangle (4.5,2.5);
\node at (-1.3,1) {$\G$};
\node at (6.3,1) {$\f:\X\to \X$};
\node[N] (1) at (1,2){\scriptsize$6$};
\node[N] (2) at (3,2){\scriptsize$1$};
\node[N] (3) at (4,1){\scriptsize$2$};
\node[N] (4) at (3,0){\scriptsize$3$};
\node[N] (5) at (1,0){\scriptsize$4$};
\node[N] (6) at (0,1){\scriptsize$5$};
\path[->]
(1) edge[red] (5)
(5) edge[red] (6)
(5) edge[very thick,Green,bend left=30] (6)
(6) edge[Green,bend left=10] (1)
(1) edge[very thick,Green, bend left=10] (6)
(4) edge[red,bend left=10] (2)
(2) edge[very thick,red] (3)
(2) edge[very thick,red,bend left=10] (4)
(3) edge[Green,bend left=10] (4)
(4) edge[Green,bend left=10] (3)
(2) edge[very thick,Green] (1)
(4) edge[very thick,Green] (1)
;
\draw[rounded corners=5,->] 
(7.8,{1-0.5}) -- (8.6,{1-1.5}) -- (8.6,{-3+1.5}) -- (7.8,{-3+0.5})
;
\node at (10.3,-1) {[~Lemma~\ref{lem:extension}~]};
\end{scope}
\begin{scope}[shift={(0,-20)}]
\draw[dotted,rounded corners=5] (-0.5,-0.5) rectangle (4.5,2.5);
\node at (-1.3,1) {$G$};
\node at (6.3,1) {$f:X\to X$};
\node[N] (1) at (1,2){\scriptsize$6$};
\node[N] (2) at (3,2){\scriptsize$1$};
\node[N] (3) at (4,1){\scriptsize$2$};
\node[N] (4) at (3,0){\scriptsize$3$};
\node[N] (5) at (1,0){\scriptsize$4$};
\node[N] (6) at (0,1){\scriptsize$5$};
\path[->]
(1) edge[red] (5)
(5) edge[red] (6)
(5) edge[Green,bend left=30] (6)
(6) edge[Green,bend left=10] (1)
(1) edge[Green, bend left=10] (6)
(4) edge[red,bend left=10] (2)
(2) edge[red] (3)
(2) edge[red,bend left=10] (4)
(3) edge[Green,bend left=10] (4)
(4) edge[Green,bend left=10] (3)
(2) edge[Green] (1)
(4) edge[Green] (1)
(5) edge[very thick,red](2)
(5) edge[very thick,Green](4)
;
\end{scope}
\end{tikzpicture}
\]
{\caption{\label{fig:case1}An illustration for the first case in the proof of Theorem~\ref{thm:main2}.}}
\end{figure}

\smallskip
\BF{Case 1: $I$ has the property $P$.}  Let $\H$ be the signed digraph obtained from $H$ by adding all the arcs $(i,j,s)$ of $G[I]$ such that $G$ has no arc from $i$ to $V\setminus I$. Let $Q$ be the subgraph of $\H$ induced by $I$ (see Figure~\ref{fig:case1}). Suppose that some connected component of $Q$ is a cycle $C$. Then $G$ has no arc from $C$ to $V\setminus I$, thus $C$ is a connected component of $G[I]$ without leaving arc. Since no connected component of $G$ is a cycle contained in $G[I]$, we deduce that $C$ has at least one entering arc. Thus $C$ is strongly connected, has no leaving arcs and at least one entering arc. This is not possible since $I$ has the property $P$. Thus no connected component of $Q$ is a cycle. Hence, by applying Theorem~\ref{thm:nilpotent2} on each connected component of $Q$, we deduce (cf. Remark~\ref{rm:component}) that there exists a degree-bounded system $g:Z\to Z$ on $Q$ such that 
\[
g^{\lambda(Q)+\beta(Q)}=\cst=(\xi_1,\dots,\xi_m).
\]

\smallskip
Let $\Y=Z_1\times \cdots\times Z_m\times Y_{m+1}\times \cdots\times Y_n$. Since $H$ has no arc from $I$ to $V\setminus I$, if $i\not\in I$ then $h_i$ only depends on variables $x_j$ with $j\not\in I$, and we can define $\h:\Y\to\Y$ as follows, without possible ambiguity:
\[
\h_i(x)=
\begin{cases}
h_i(x_{m+1},\dots,x_n)&\textrm{ if $i\not\in I$}\\
g_i(x_1,\dots,x_m)&\textrm{ if $i\in I$.}
\end{cases}
\]

\smallskip
By construction, $\h$ is a degree-bounded system on $\H$ such that: 
\begin{enumerate}
\item[(1)]
$\h_i(\Y)\subseteq h_i(Y)$ for all $i\not\in I$,
\item[(2)]
$\h_i^{\lambda(Q)+\beta(Q)}(\Y)\subseteq \{\xi_i\}=Y_i=h_i(Y)$ for all $i\in I$,
\item[(3)]
$\h(x)=h(x)$ for all $x\in Y$.
\end{enumerate}

\smallskip
Let $\G$ be the signed digraph obtained from $\H$ by adding all the arcs $(i,j,s)$ of $G$ with $j\not\in I$ (see Figure~\ref{fig:case1}). Let $A$ be the set of sources of $\H$ that are not sources of $\G$, and let us prove that $A=\emptyset$. Suppose that $j$ is a source of $\H$. If $j\in I$ then $j$ is obviously a source of $\G$, and if $j\notin I$ then $j$ is a source of $H\setminus I$. Hence, by hypothesis, $j$ is a source of $G$ and thus of $\G$ too. Thus, $A=\emptyset$. This implies that $\G$ has no arc from a sink of $\H$ to a source of $\H$. Furthermore, $\G$ has no arc from a vertex of $\H$ to an isolated vertex $i$ of $\H$ since otherwise $i\not\in I$, thus $i$ is isolated in $H$ but not in $G$, a contradiction. Thus $\G$ and $\H$ satisfy the conditions of Lemma~\ref{lem:extension}. Let $B$ be the set of sinks of $\H$ that are not sinks of $\G$. If $i$ is a sink of $\H$ and $i\not\in I$, then $i$ is a sink of $H\setminus I$ and, by hypothesis, $i$ is sink of $G$ and thus of $\G$ too. Thus $B\subseteq I$. Let $i\in I$ and $j\in B$. Then $G$ has an arc from $j$ to $V\setminus I$, hence, by definition, $Q$ has no arc from $j$ to $i$ and we deduce that $\G$ has no arc from $j$ to $i$. Thus $\G_i\cap B=\emptyset$ for every $i\in I$. We then deduce from Lemma~\ref{lem:extension}, that there exists a degree-bounded system $\f:\X\to\X$ on $\G$ such that      
\begin{enumerate}
\item[(4)] 
$\f(\X)\subseteq \h(\Y)$,    
\item[(5)]
$\f(x)=\h(x)$ for all $x\in\Y$ such that $x_i=\xi_i$ for all $i\in I$, that is, for all $x\in Y$,  
\item[(6)]
$\f_i(x)=\h_i(x)$ for all $x\in\Y$ and $i\in I$.
\end{enumerate}

\smallskip
As a consequence, $\f$ has the following four properties:
\begin{enumerate}
\item[(7)]
$\f_i(\X)\subseteq h_i(Y)$ for all $i\not\in I$, 
\item[(8)]
$\f^\ell_i(\X)\subseteq \h^\ell_i(\Y)$ for all $i\in I$ and $\ell\geq 1$,
\item[(9)]
$\f^{\lambda(Q)+\beta(Q)}_i(\X)\subseteq h_i(Y)$ for all $i\in I$,
\item[(10)]
$\f(x)=h(x)$ for all $x\in Y$. 
\end{enumerate}
Indeed, property (7) results from (1) and (4), and property (10) results from (5) and (3). We now prove (8) by induction on $\ell$. The case $\ell=1$ is given by (4), and if $\ell\geq 2$ and $i\in I$ then 
\[
\f^\ell_i(\X)=\f_i(\f^{\ell-1}(\X))\subseteq \f_i(\h^{\ell-1}(\Y))=\h^\ell_i(\Y),
\]
where the inclusion results from the induction hypothesis and the fact that $\f_i$ only depends on variables with indices in $I$ (because $\G$ has no arc from $V\setminus I$ to $I$), and where the last equality results from (6). Finally, (9) results from (8) and (2). 

\smallskip
Let $A$ be the set of sources of $\G$ that are not sources of $G$, and let $B$ be the set of sinks of $\G$ that are not sinks of $G$. Since $G$ is obtained from $\G$ by adding only arcs $(i,j,s)$ with $j\in I$ (see Figure~\ref{fig:case1}), we have $A\subseteq I$. Let us prove that $B$ is empty. Let $i$ be a sink of $\G$. If $i\not\in I$ then $i$ is a sink of $H\setminus I$ and thus, by hypothesis, $i$ is a sink of $G$. If $i\in I$ then, by the definition of $\G$, $G$ has no arc from $i$ to $V\setminus I$ thus $i$ is a sink of $G$, and since $i$ is a sink of $\H$, by the definition of $\H$, $G$ has no arc from $i$ to $I$. Thus $i$ is a sink of $G$. This proves that $B=\emptyset$. We deduce that $G$ has no arc from a sink of $\G$ to a source of $\G$. Now, suppose, for a contradiction, that $G$ has an arc from a vertex of $\G$ to an isolated vertex $i$ of $\G$. If $i\not\in I$ then $i$ is isolated in $H\setminus I$ but not in $G$, a contradiction, thus $i\in I$. Hence, $i$ is a trivial connected component of $G[I]$, which is strongly connected, which has no leaving arcs, and which has at least one entering arc, and this contradicts the fact that $I$ has the property $P$. Thus $G$ has no arc from a vertex of $\G$ to an isolated vertex of $\G$. Thus $G$ and $\G$ satisfy the conditions of Lemma~\ref{lem:extension}, with $A\subseteq I$ and $B=\emptyset$. Therefore, there exists a degree-bounded system $f:X\to X$ on $G$ such that 
\begin{enumerate}
\item[(11)]
$f(X)\subseteq \X$,
\item[(12)]
$f_i(X)\subseteq \f_i(\X)$ for all $i\not \in A$,
\item[(13)]
$f(x)=\f(x)$ for all $x\in\X$.
\end{enumerate}

\smallskip
As a consequence, $f$ has the following properties:
\begin{itemize}
\item[(14)] 
$f_i(X)\subseteq h_i(Y)$ for all $i\not\in I$, 
\item[(15)] 
$f^{|I|+1}_i(X)\subseteq h_i(Y)$ for all $i\in I$,
\item[(16)] 
$f(x)=h(x)$ for all $x\in Y$.
\end{itemize}
Indeed, (14) results from (12), (7) and the fact that $A\subseteq I$. Also, (16) results from (13) and (10). It remains to prove (15). Suppose first that $A=\emptyset$. Then following (12) we have  $f(X)\subseteq \f(\X)$ and we deduce from (13) that $f^\ell(X)\subseteq \f^\ell(\X)$ for all $\ell\geq 1$. Thus by (9) we have 
\[
f^{\lambda(Q)+\beta(Q)}_i(X)\subseteq h_i(Y). 
\]
Since $\lambda(Q)\leq|I|$ and $\beta(Q)\leq 1$ we obtain (15). Suppose now that $A\neq\emptyset$. Then following (11) and (13) we have $f^{\ell+1}(X)\subseteq \f^\ell(\X)$ for all $\ell\geq 0$. Thus by (9) we have  
\[
f^{\lambda(Q)+\beta(Q)+1}_i(X)\subseteq h_i(Y).
\]
If $Q$ is basic then $\beta(Q)=0$ and (15) follows. Otherwise, $Q$ is not basic, and since $A\neq\emptyset$ and $A\subseteq I$, $Q$ has a source. Hence, $Q$ has at least two distinct initial strong components. Thus, as explained in Remark~\ref{rem:lambda}, we have $\lambda(Q)<|I|$ and (15) follows.  

\smallskip
We deduce from (14), (15) and (16) that $f$ converges toward $h$ in at most $|I|+1$ steps.

\begin{figure}
\[
\begin{tikzpicture}
\draw[dotted,rounded corners=5] (-0.75,-0.75) rectangle (3.75,2.25);
\node at (-1.55,0.75) {$G$};
\node[N] (4) at (0.0,0.0){\scriptsize$4$};
\node[N] (5) at (1.5,0.0){\scriptsize$5$};
\node[N] (6) at (3.0,0.0){\scriptsize$6$};
\node[N] (1) at (0.0,1.5){\scriptsize$1$};
\node[N] (2) at (1.5,1.5){\scriptsize$2$};
\node[N] (3) at (3.0,1.5){\scriptsize$3$};
\path[->,thick]
(1) edge[Green,bend left=20] (4)
(4) edge[Green,bend left=20] (1)
(2) edge [red] (1)
(4) edge [Green] (5)
(2) edge[Green,bend left=20] (3)
(3) edge[Green,bend left=20] (2)
(5) edge[Green,bend left=20] (6)
(6) edge[red,bend left=20] (5)
;
\begin{scope}[shift={(0,-4)}]
\node at (-1.55,0.75) {$H$};
\node at (5.5,0.75) {$h:Y\to Y$};
\draw[dotted,rounded corners=5] (-0.75,-0.75) rectangle (3.75,2.25);
\draw[dotted,rounded corners=5] (1.05,1.05) rectangle (3.45,1.95);
\node[preaction={fill,white}] at (2.25,1.05) {\scriptsize $\I$};
\draw[dotted,rounded corners=5] (1.05,-0.45) rectangle (3.45,0.45);
\node[preaction={fill,white}] at (2.25,0.45) {\scriptsize $J$};
\draw[dotted,rounded corners=5] (0.9,-0.6) rectangle (3.6,2.1);
\node[preaction={fill,white}] at (0.9,0.75) {\small $I$};
\node[N] (4) at (0.0,0.0){\scriptsize$4$};
\node[N] (5) at (1.5,0.0){\scriptsize$5$};
\node[N] (6) at (3.0,0.0){\scriptsize$6$};
\node[N] (1) at (0.0,1.5){\scriptsize$1$};
\node[N] (2) at (1.5,1.5){\scriptsize$2$};
\node[N] (3) at (3.0,1.5){\scriptsize$3$};
\path[->,thick]
(1) edge[Green,bend left=20] (4)
(4) edge[Green,bend left=20] (1)
;
\draw[rounded corners=5,->] 
(7,{0.75-0.5}) -- (7.8,{0.75-1.5}) -- (7.8,{-3.3+1.5}) -- (7,{-3.3+0.5})
;
\node at (10.3,{-2+0.75}) {[~first case of Theorem~\ref{thm:main2}~]};
\end{scope}
\begin{scope}[shift={(0,-8)}]
\node at (-1.55,0.75) {$\G$};
\node at (5.5,0.75) {$\f:\X\to\X$};
\draw[dotted,rounded corners=5] (-0.75,-0.75) rectangle (3.75,2.25);
\node[N] (4) at (0.0,0.0){\scriptsize$4$};
\node[N] (5) at (1.5,0.0){\scriptsize$5$};
\node[N] (6) at (3.0,0.0){\scriptsize$6$};
\node[N] (1) at (0.0,1.5){\scriptsize$1$};
\node[N] (2) at (1.5,1.5){\scriptsize$2$};
\node[N] (3) at (3.0,1.5){\scriptsize$3$};
\path[->,thick]
(1) edge[Green,bend left=20] (4)
(4) edge[Green,bend left=20] (1)
(2) edge [red] (1)
(2) edge[Green,bend left=20] (3)
(3) edge[Green,bend left=20] (2)
;
\draw[rounded corners=5,->] 
(7,{0.75-0.5}) -- (7.8,{0.75-1.5}) -- (7.8,{-7.3+1.5}) -- (7,{-7.3+0.5})
;
\end{scope}
\begin{scope}[shift={(0,-12)}]
\draw[dotted,rounded corners=5] (-0.75,-0.75) rectangle (3.75,2.25);
\node at (-1.55,0.75) {$Q$};
\node at (5.5,0.75) {$g:Z\to Z$};
\node at (5.5,0) {[~Theorem~\ref{thm:nilpotent2}~]};
\node[N] (4) at (0.0,0.0){\scriptsize$4$};
\node[N] (5) at (1.5,0.0){\scriptsize$5$};
\node[N] (6) at (3.0,0.0){\scriptsize$6$};
\node[N] (1) at (0.0,1.5){\scriptsize$1$};
\node[N] (2) at (1.5,1.5){\scriptsize$2$};
\node[N] (3) at (3.0,1.5){\scriptsize$3$};
\path[->,thick]
(4) edge [Green] (5)
(5) edge[Green,bend left=20] (6)
(6) edge[red,bend left=20] (5)
;
\draw[rounded corners=5,-] 
(7,{0.75-0.5}) -- (7.8,{0.75-1.5}) -- (7.8,-1)
;
\end{scope}
\begin{scope}[shift={(0,-16)}]
\draw[dotted,rounded corners=5] (-0.75,-0.75) rectangle (3.75,2.25);
\node at (-1.55,0.75) {$G$};
\node at (5.5,0.75) {$f:X\to X$};
\node[N] (4) at (0.0,0.0){\scriptsize$4$};
\node[N] (5) at (1.5,0.0){\scriptsize$5$};
\node[N] (6) at (3.0,0.0){\scriptsize$6$};
\node[N] (1) at (0.0,1.5){\scriptsize$1$};
\node[N] (2) at (1.5,1.5){\scriptsize$2$};
\node[N] (3) at (3.0,1.5){\scriptsize$3$};
\path[->,thick]
(1) edge[Green,bend left=20] (4)
(4) edge[Green,bend left=20] (1)
(2) edge [red] (1)
(4) edge [Green] (5)
(2) edge[Green,bend left=20] (3)
(3) edge[Green,bend left=20] (2)
(5) edge[Green,bend left=20] (6)
(6) edge[red,bend left=20] (5)
;
\end{scope}
\end{tikzpicture}
\]
{\caption{\label{fig:case2}An illustration for the second case in the proof of Theorem~\ref{thm:main2}.}}
\end{figure}
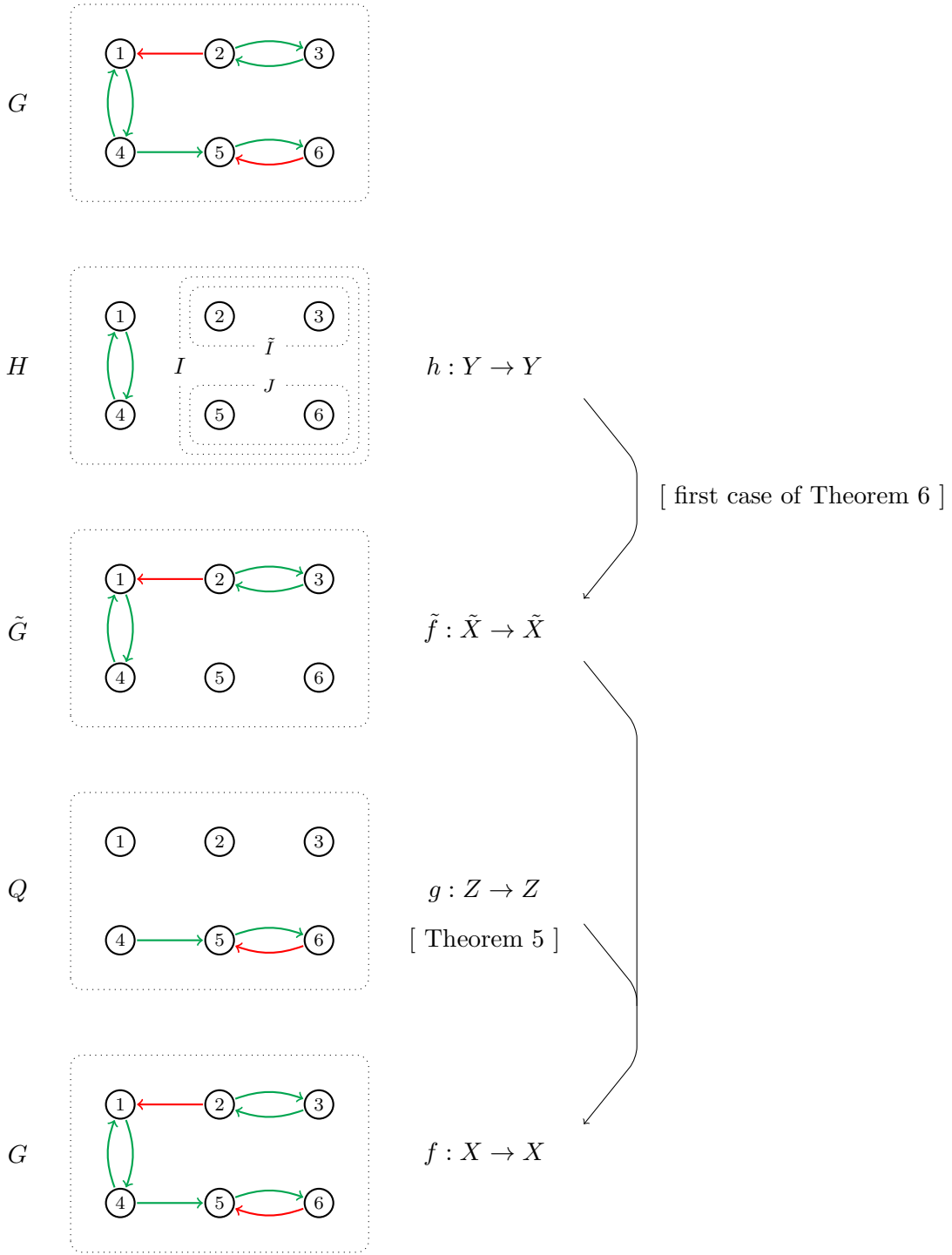

\smallskip
\BF{Case 2: $I$ has not the property $P$.} 
Let $I_1$ be the set of vertices that belongs to a connected component $F$ of $G[I]$ such that $F$ is not strongly connected, or has a leaving arc, or has no entering arcs. Let $J=I\setminus \I$. Since $I$ has not the property $P$, $J$ is not empty. Let $\G$ be the signed digraph obtained from $G$ by removing every arc $(i,j,s)$ with $j\in J$ (see Figure \ref{fig:case2}). Then, $\I$ is the set of vertices isolated in $H$ but not in $\G$. Furthermore, it is straightforward to check that no connected component of $\G$ is a cycle contained in $\G[\I]$ and that every source (resp. sink) of $\G\setminus\I$ is a source (resp. sink) of $\G$. Furthermore, since $\I$ has the property $P$ with respect to $\G$, according to the first case, there exists a degree-bounded system $\f:\X\to\X$ on $\G$ such that
\begin{itemize}
\item[(1)] 
$\f_i(\X)\subseteq h_i(Y)$ for all $i\not\in \I$,
\item[(2)] 
$\f^{|\I|+1}_i(\X)\subseteq h_i(Y)$ for all $i\in \I$, and thus $\f^{|I|+1}_i(\X)\subseteq h_i(Y)$ for all $i\in \I$,
\item[(3)] 
$\f(x)= h(x)$ for all $x\in Y$.
\end{itemize}

\smallskip 
Let $\xi\in\Z^n$ be defined by $\xi_i=\max(\X_i)$ for all $i\in \{1,\dots,n\}$. Let $J=I\setminus\I$ and let $Q$ be the spanning subgraph of $G$ that only contains the arcs $(i,j,s)$ of $G$ with $j\in J$. So $\G$ and $Q$ are arc-disjoint and $G=\G\cup Q$ (see Figure~\ref{fig:case2}). Since every connected components $G[J]$ is strongly connected, has no leaving arcs and has an entering arc, we deduce that every connected component of $Q$ is basic. Hence, by applying Theorem~\ref{thm:nilpotent2} on each connected component of $H$ (cf. Remark~\ref{rm:component}) we deduce that there exists a degree-bounded system $g:Z\to Z$ on $Q$ such that 
\[
g^{\lambda(Q)}=\cst=\xi.
\]
and such that $\xi_i=\min(Z_i)$ for all the sources $i$ of $Q$. Since if $i\not\in J$ then $i$ is a source of $Q$ we have $\max(\X_i)=\xi_i=\min(Z_i)$ for all $i\not\in J$. Remark also that $Y_i=\X_i=\{\xi_i\}\subseteq Z_i$ for all $i\in J$. 

\smallskip
Let $X$ be defined by:
\[
X=X_1\times\dots\times X_n
\quad\text{with}\quad
X_i=
\begin{cases}
\X_i\cup Z_i&\text{if }i\not\in J\\
Z_i&\text{if }i\in J.
\end{cases}
\]
Let $\pi:X\to Z$ and $\tpi:X\to \X$ be the surjections defined by: 
\[
\pi_i(x)=
\begin{cases}
\max(x_i,\xi_i)&\text{if }i\not\in J\\
x_i&\text{if }i\in J.
\end{cases}
\quad\text{and}\quad
\tpi_i(x)=
\begin{cases}
\min(x_i,\xi_i)&\text{if }i\not\in J\\
\xi_i&\text{if }i\in J
\end{cases}
\]
Finally, let $f:X\to X$ be defined by:
\[
f_i(x)=
\begin{cases}
\f_i(\tpi(x))&\textrm{if }i\notin J\\
g_i(\pi(x))&\textrm{if }i\in J.
\end{cases}
\]
Since $\f$ is a degree-bounded system on $\G$ and $g$ is a degree-bounded system on $Q$, it is easy to check that $f$ is a degree-bounded system on $G=\G\cup Q$. $f$ has also the following properties:
\begin{itemize}
\item[(4)]
$f^\ell_i(X)\subseteq\f^\ell_i(\X)$ for all $i\not\in J$ and $\ell\geq 1$.  
\item[(5)]
$f^\ell_i(x)=g^\ell_i(\pi(x))$ for all $i\in J$, $x\in X$ and $\ell\geq 1$, so that  
$f^{\lambda(Q)}_i=\cst=\xi_i$ for all $i\in J$,
\item[(6)]
$f(x)=\f(x)$ for all $x\in\X$.
\end{itemize}
We first prove (5) by induction on $\ell$. The case $\ell=1$ results from the definition of $f_i$ for $i\in J$. If $\ell\geq 2$ then, by definition, we have 
\[\tag{$*$}\label{eq:h}
f^{\ell}_i(x)=f_i(f^{\ell-1}(x))=g_i(\pi(f^{\ell-1}(x))).
\] 
By induction hypothesis, for all $j\in J$ we have $f^{\ell-1}_j(x)=g^{\ell-1}_j(\pi(x))\geq\xi_j$ thus $\pi_j(f^{\ell-1}(x))=g^{\ell-1}_j(\pi(x))$. If $j\not\in J$ then $f^{\ell-1}_j(x)\leq\xi_j$ thus $\pi_j(f^{\ell-1}(x))=\xi_j$. Since $j$ is a source of $Q$ and $\xi$ is a fixed point of~$g$, we have $g_j=\cst=\xi_j$. Thus, $\pi_j(f^{\ell-1}(x))=g^{\ell-1}_j(\pi(x))$. Consequently, $\pi(f^{\ell-1}(x))=g^{\ell-1}(\pi(x))$ and we deduce from (\ref{eq:h}) that $f^\ell_i(x)=g^\ell_i(\pi(x))$, completing the induction step. We now prove (6). Let $x\in\X$. Then $\tpi(x)=x$ thus $f_i(x)=\f_i(x)$ for all $i\notin J$. Furthermore, $\pi(x)=\xi$ and if $i\in I$ then $\X_i=\{\xi_i\}=\{x_i\}$. Since $\xi$ is a fixed point of $g$ we deduce that $f_i(x)=g_i(\xi)=\xi_i=\f_i(x)$. Hence $f(x)=\f(x)$ and (6) is proved. We finally prove (4) by induction on $\ell$. The case $\ell=1$ is an immediate consequence of the definition of $f_i$ for $i\not\in J$. If $\ell\geq 2$ and $i\not\in J$ then 
\[
f^\ell_i(X)=f_i(\f^{\ell-1}(\X))\subseteq f_i(\f^{\ell-1}(\X))=\f^\ell_i(\X),
\]
where the inclusion results from the induction hypothesis and the fact that $f_i$ only depends on variables with indices which are not in $J$ (that is, $\G$ has no arc from $J$ to $V\setminus J$), and where the last equality results from (6).  

\smallskip
As a consequence, $f$ has the following three properties:
\begin{itemize}
\item[(7)] $f_i(X)\subseteq h_i(Y)$ for all $i\not\in I$. 
\item[(8)] $f^{|I|+1}_i(X)\subseteq h_i(Y)$ for all $i\in I$.
\item[(9)] $f(x)=h(x)$ for all $x\in Y$.
\end{itemize}
Indeed, (7) results from (4) and (1), and (9) results from (6) and (3). It remains to prove (8). From $(4)$ and $(2)$ we have $f^{|I|+1}_i(X)\subseteq h_i(Y)$ for all $i\in\I$, and from (5) we have $f^\ell_i(X)=\{\xi_i\}=h_i(Y)$ for all $i\in J$ and $\ell\geq \lambda(Q)$. We deduce that is sufficient to prove that $\lambda(Q)\leq |I|+1$. Since $J$ corresponds to the set of non-sources of $Q$ and $Q$ is basic, by Remark~\ref{rem:lambda}, we have $\lambda(Q)\leq |J|+1\leq |I|+1$. This completes the proof of (8). 

\smallskip
We deduce from (7), (8) and (9) that $f$ converges toward $h$ in at most $|I|+1$ steps.
\end{proof}

We now prove Lemma~\ref{lem:extension}, using the following notations. For all $X\subseteq \mathbb{Z}^n$ and $y\in\mathbb{Z}^n$, $X+y=\{x+y\,|\,x\in X\}$. We denote by $e_i$ the $n$-tuple whose components are all equal to $0$, excepted th $i$th one, which is equal to $1$.

\begin{proof}[{\bf Proof of Lemma~\ref{lem:extension}}]
We proceed by induction on the number of arcs in $G$ that are not in $H$. If no such arcs exist, then $G=H$ and there is nothing to prove. So assume that $G$ has an arc 
\[
a=(j,i,s)
\]
that is not in $H$. Let $\G=G\setminus a$ be the spanning subgraph obtained from $G$ by removing $a$. Let $\A$ be the set of sources of $H$ that are not sources of $\G$ and let $\B$ be the set of sinks of $H$ that are no sinks of $\G$. By induction hypothesis, there exists a degree-bounded system $\f:\X\to\X$ on $\G$ such that 
\begin{enumerate}
\item[(i)] 
$\f(\X)\subseteq Y$,
\item[(ii)] 
$\f_k(\X)\subseteq h_k(Y)$ for all $k\not\in \A$.
\item[(iii)]
$\f(x)=h(x)$ for all $x\in Y$ such that $x_k=\xi_k$ for all $k\in \B$.   
\item[(iv)]
$\f_k(x)=h_k(x)$ for all $x\in Y$ and vertex $k$ such that $\G_k\cap \B=\emptyset$.
\end{enumerate}
Since $G$ has no arc from a source of $H$ to a sink of $H$, $G$ has no arc from a source of $\G$ to a sink of $\G$. Furthermore, since $G$ has no arc from a vertex of $H$ to an isolated vertex of $H$, $G$ has no arc from a vertex of $\G$ to an isolated vertex of $\G$. This leaves the following three possibilities.

\smallskip
\BF{Case 1: $j$ is not a sink of $\G$ and $i$ is not a source of $\G$.} 
\begin{itemize}
\item
If $s=+$ we set $X=\X\cup (\X+e_j)$ and define $f:X\to X$ as follows: 
\[
\begin{array}{lcl}
f_i(x)&=&
\left\{
\begin{array}{ll}
\max(\f_i(\X))&\textrm{if $x_j=\max(X_j)$}\\
\f_i(x)&\textrm{otherwise}\\
\end{array}
\right.
\\[5mm]
f_k(x)&=&
\left\{
\begin{array}{ll}
\f_k(x-e_j)&\textrm{if $x\not\in \X$}\\
\f_k(x)&\textrm{otherwise}\\
\end{array}
\right.
\quad\forall k\neq i.
\end{array}
\]
Hence, $f(X)\subseteq\f(\X)$ and $f(x)=\f(x)$ for all $x\in \X$. Since $i$ is a not a sink of $\G$ we have $A=\A$ and since $j$ is not a source of $\G$, we have $B=\B$. Thus, using the properties (i)-(iv), we deduce that $f$ satisfies the four points of the statement. It remains to prove that $f$ is a degree-bounded system on $G$. 

\smallskip
Let $G_f$ be the interaction graph of $f$. Since $f(x)=\f(x)$ for all $x\in \X$, $\G$ is a subgraph of $G_f$. Since it is straightforward to prove that $G_f$ is a subgraph of $G$, to prove that $G_f=G$, it is sufficient to prove that $a$ is an arc of $G_f$. Let $x\in X$ be such that $f_i(x)=\min(\f_i(\X))$ and suppose that $x_j$ is maximal for this property. Since $i$ is not a source of $\G$, we have $\min(\f_i(\X))<\max(\f_i(\X))$ and we deduce from the definition of $f_i$ that $x_j<\max(X_j)$. Consequently, $f_i(x+e_j)>\min(\f_i(\X))$ so $G_f$ has a positive arc from $j$~to~$i$. Thus $G_f=G$, that is, $f$ is an FDS on $G$. Since $\f$ is degree-bounded, it is straightforward to prove that $f$ is degree-bounded too. 

\item
If $s=-$, we proceed similarly with $X=\X\cup (\X+e_j)$ and $f:X\to X$ defined by: 
\[
\begin{array}{lcl}
f_i(x)&=&
\left\{
\begin{array}{ll}
\min(\f_i(\X))&\textrm{if $x_j=\max(X_j)$}\\
\f_i(x)&\textrm{otherwise}\\
\end{array}
\right.
\\[5mm]
f_k(x)&=&
\left\{
\begin{array}{ll}
\f_k(x-e_j)&\textrm{if $x\not\in \X$}\\
\f_k(x)&\textrm{otherwise}\\
\end{array}
\right.
\quad\forall k\neq i.
\end{array}
\]
\end{itemize}

\smallskip
\BF{Case 2: $j$ is a sink of $\G$ and $i$ is not a source of $\G$.} 
Since $j$ is a sink of $\G$ we have $|\X_j|=1$ if $j$ is isolated in $\G$, and $|\X_j|=2$ otherwise. Let 
\[
X=
\begin{cases}
\X&\text{if }|\X_j|= 2\\
\X\cup(\X+e_j)&\text{if }|\X_j|= 1\\
\end{cases}
\]
Hence $|X_j|=2$, so either $\xi_j<\max(X_j)$ or $\xi_j>\min(X_j)$, and we consider four cases.
\begin{itemize}
\item
If $s=+$ and $\xi_j<\max(X_j)$, we define $f:X\to X$ by:
\[
\begin{array}{lcl}
f_i(x)&=&
\left\{
\begin{array}{ll}
\max(\f_i(\X))&\textrm{if $x_j=\max(X_j)$}\\
\f_i(x)&\textrm{otherwise}\\
\end{array}
\right.
\\[5mm]
f_k(x)&=&
\left\{
\begin{array}{ll}
\f_k(x-e_j)&\textrm{if $x\not\in \X$}\\
\f_k(x)&\textrm{otherwise}\\
\end{array}
\right.
\quad\forall k\neq i.
\end{array}
\]
Hence, $f(X)\subseteq\f(\X)$ and $f(x)=\f(x)$ for all $x\in X$ with $x_j=\xi_j$. Since $i$ is a not a source of $\G$ we have $A=\A$, we deduce from (i) and (ii) that $f$ satisfies the first two points of the statement. Since $B=\B\cup\{j\}$ and $j\in G_i\cap B$, we then deduce from (iii) and (iv) that $f$ satisfies the last two points. It remains to prove that $f$ is a degree-bounded system on $G$. We prove that $G$ is the interaction graph of $f$ exactly as in the previous case. Next, since $\f$ is a degree-bounded system on $\G$, and since $|X_j|=2=\dout_G(j)+1$, we deduce that $f$ is degree-bounded. Thus $f$ is indeed a degree-bounded system on $G$.  

\item
If $s=+$ and $\xi_j>\min(X_j)$, we proceed similarly with $f:X\to X$ defined by:
\[
\begin{array}{lcl}
f_i(x)&=&
\left\{
\begin{array}{ll}
\min(\f_i(\X))&\textrm{if $x_j=\min(X_j)$}\\
\f_i(x)&\textrm{otherwise}\\
\end{array}
\right.
\\[5mm]
f_k(x)&=&
\left\{
\begin{array}{ll}
\f_k(x-e_j)&\textrm{if $x\not\in \X$}\\
\f_k(x)&\textrm{otherwise}\\
\end{array}
\right.
\quad\forall k\neq i.

\end{array}
\]
\item
If $s=-$ and $\xi_j<\max(X_j)$, we proceed similarly with $f:X\to X$ defined by:
\[
\begin{array}{lcl}
f_i(x)&=&
\left\{
\begin{array}{ll}
\min(\f_i(\X))&\textrm{if $x_j=\max(X_j)$}\\
\f_i(x)&\textrm{otherwise}\\
\end{array}
\right.
\\[5mm]
f_k(x)&=&
\left\{
\begin{array}{ll}
\f_k(x-e_j)&\textrm{if $x\not\in \X$}\\
\f_k(x)&\textrm{otherwise}\\
\end{array}
\right.
\quad\forall k\neq i.

\end{array}
\]
\item
If $s=-$ and $\xi_j>\min(X_j)$, we proceed similarly with $f:X\to X$ defined by:
\[
\begin{array}{lcl}
f_i(x)&=&
\left\{
\begin{array}{ll}
\max(\f_i(\X))&\textrm{if $x_j=\min(X_j)$}\\
\f_i(x)&\textrm{otherwise}\\
\end{array}
\right.
\\[5mm]
f_k(x)&=&
\left\{
\begin{array}{ll}
\f_k(x-e_j)&\textrm{if $x\not\in \X$}\\
\f_k(x)&\textrm{otherwise}\\
\end{array}
\right.
\quad\forall k\neq i.
\end{array}
\]
\end{itemize}

\smallskip
\BF{Case 3: $j$ is not a sink of $\G$ and $i$ is a non-isolated source of $\G$.}
Then $i$ is a source of both $\G$ and $H$, thus $\f_i=\mathrm{cst}=c_i$ and $h_i=\cst=c'_i$, and following (iii) $c_i=c'_i$. Since $i$ is not isolated in $\G$ we have $|\X_i|\geq 2$ and thus either $c_i<\max(\X_i)$ or $c_i>\min(\X_i)$, and we consider the following four cases. In every cases we set 
\[
X=\X\cup (\X+ e_j).
\]

\begin{itemize}
\item
If $s=+$ and $c_i< \max(\X_i)$, we define $f:X\to X$ by:
\[
\begin{array}{lcl}
f_i(x)&=&
\left\{
\begin{array}{ll}
\max(\X_i)&\textrm{if $x_j=\max(X_j)$}\\
c_i&\textrm{otherwise}\\
\end{array}
\right.
\\[5mm]
f_k(x)&=&
\left\{
\begin{array}{ll}
\f_k(x-e_j)&\textrm{if $x\not\in\X$}\\
\f_k(x)&\textrm{otherwise}\\
\end{array}
\right.
\quad\forall k\neq i.
\end{array}
\]

Hence, $f_i(X)\subseteq \X_i$ and $f_k(X)\subseteq \f_k(\X)$ for all $k\neq i$. Since $A=\A\cup\{i\}$ we deduce from (i) and (ii) that $f$ satisfies the first two points of the statement.  Since $B=\B$ and $h(x)=\f(x)$ for all $x\in \X$, we deduce from (iii) and (iv) that $f$ satisfies the last two points of the statement. It remains to prove that $f$ is a degree-bounded system on $G$.  

\smallskip
Let $G_f$ be the interaction graph of $f$. Since $f(x)=\f(x)$ for all $x\in \X$, $\G$ is a subgraph of $G_f$. Since it is straightforward to show that $G_f$ is a subgraph of $G$, to prove that $G_f=G$, it is sufficient to prove that $a$ is an arc of $G_f$. Let $x\in X$ with $x_j=\max(X_j)-1$. Then $f_i(x)=c_i<\max(\X_i)$ and $f_i(x+e_j)=\max(\X_i)$ thus $a$ is indeed an arc of $G_f$. Thus $f$ is a FDS on $G$, and since $\f$ is degree-bounded, it is straightforward to prove that $f$ is degree-bounded too. 

\item
If $s=+$ and $c_i>\min(\X_i)$ we proceed similarly with $f:X\to X$ define by:
\[
\begin{array}{rcl}
f_i(x)&=&
\left\{
\begin{array}{ll}
c_i&\textrm{if $x_j=\max(X_j)$}\\
\min(\X_i)&\textrm{otherwise}\\
\end{array}
\right.
\\[5mm]
f_k(x)&=&
\left\{
\begin{array}{ll}
\f_k(x-e_j)&\textrm{if $x\not\in\X$}\\
\f_k(x)&\textrm{otherwise}\\
\end{array}
\right.
\quad\forall k\neq i.
\end{array}
\]
\item
If $s=-$ and $c_i<\max(\X_i)$, we proceed similarly with $f:X\to X$ defined by:
\[
\begin{array}{rcl}
f_i(x)&=&
\left\{
\begin{array}{ll}
c_i&\textrm{if $x_j=\max(X_j)$}\\
\max(\X_i)&\textrm{otherwise}\\
\end{array}
\right.
\\[5mm]
f_k(x)&=&
\left\{
\begin{array}{ll}
\f_k(x-e_j)&\textrm{if $x\not\in\X$}\\
\f_k(x)&\textrm{otherwise}\\
\end{array}
\right.
\quad\forall k\neq i.
\end{array}
\]
\item
If $s=-$ and $c_i>\min(X_i)$, we proceed similarly with $f:X\to X$ defined by:
\[
\begin{array}{rcl}
f_i(x)&=&
\left\{
\begin{array}{ll}
\min(\X_i)&\textrm{if $x_j=\max(X_j)$}\\
c_i&\textrm{otherwise}\\
\end{array}
\right.
\\[5mm]
f_k(x)&=&
\left\{
\begin{array}{ll}
\f_k(x-e_j)&\textrm{if $x\not\in\X$}\\
\f_k(x)&\textrm{otherwise}\\
\end{array}
\right.
\quad\forall k\neq i.
\end{array}
\]
\end{itemize}
\end{proof}

\paragraph{Acknowledgment} I wish to thank Emmanuel Bonan and Tony Silva for stimulating discussions.

\bibliographystyle{plain}
\bibliography{BIB}

\begin{thebibliography}{10}

\bibitem{AMN14}
W.~Abou-Jaoud{\'e}, P.~Monteiro, A.~Naldi, M.~Grandclaudon, V.~Soumelis,
  C.~Chaouiya, and D.~Thieffry.
\newblock Model checking to assess t-helper cell plasticity.
\newblock {\em Frontiers in bioengineering and biotechnology}, 2, 2014.

\bibitem{A08}
J.~Aracena.
\newblock Maximum number of fixed points in regulatory {B}oolean networks.
\newblock {\em Bulletin of Mathematical Biology}, 70(5):1398--1409, 2008.

\bibitem{ARS17}
J.~Aracena, A.~Richard, and L.~Salinas.
\newblock Number of fixed points and disjoint cycles in monotone boolean
  networks.
\newblock {\em SIAM Journal on Discrete Mathematics}, 31(3):1702--1725, 2017.

\bibitem{CD02}
O.~Cinquin and J.~Demongeot.
\newblock Positive and negative feedback: strinking a balance between necessary
  antagonists.
\newblock {\em Journal of Theoretical Biology}, 216:229--241, 2002.

\bibitem{GRF16}
M.~Gadouleau, A.~Richard, and E.~Fanchon.
\newblock Reduction and fixed points of boolean networks and linear network
  coding solvability.
\newblock {\em IEEE Transactions on Information Theory}, 62(5):2504--2519,
  2016.

\bibitem{GR11}
M.~Gadouleau and S.~Riis.
\newblock Graph-theoretical constructions for graph entropy and network coding
  based communications.
\newblock {\em IEEE Transactions on Information Theory}, 57(10):6703--6717,
  2011.

\bibitem{GR16}
Maximilien Gadouleau and Adrien Richard.
\newblock Simple dynamics on graphs.
\newblock {\em Theoretical Computer Science}, 628:62--77, 2016.

\bibitem{GM90}
E.~Goles and S.~Mart{\'\i}nez.
\newblock {\em Neural and Automata Networks: Dynamical Behavior and
  Applications}.
\newblock Kluwer Academic Publishers, 1990.

\bibitem{G98}
J.L. Gouz\'e.
\newblock Positive and negative circuits in dynamical systems.
\newblock {\em Journal of Biological Systems}, 6:11--15, 1998.

\bibitem{H82}
J.~Hopfield.
\newblock Neural networks and physical systems with emergent collective
  computational abilities.
\newblock {\em Proc. Nat. Acad. Sc. U.S.A.}, 79:2554--2558, 1982.

\bibitem{K69}
S.~A. Kauffman.
\newblock Metabolic stability and epigenesis in randomly connected nets.
\newblock {\em Journal of Theoretical Biology}, 22:437--467, 1969.

\bibitem{K93}
S.~A. Kauffman.
\newblock {\em Origins of Order Self-Organization and Selection in Evolution}.
\newblock Oxford University Press, 1993.

\bibitem{KS19}
M.~Kaufman and C.~Soul{\'e}.
\newblock On the multistationarity of chemical reaction networks.
\newblock {\em Journal of theoretical biology}, 465:126--133, 2019.

\bibitem{KST07}
M.~Kaufman, C.~Soul\'e, and R.~Thomas.
\newblock A new necessary condition on interaction graphs for
  multistationarity.
\newblock {\em Journal of Theoretical Biology}, 248(4):675--685, 2007.

\bibitem{N15}
N.~Le~Nov{\`e}re.
\newblock Quantitative and logic modelling of molecular and gene networks.
\newblock {\em Nature Reviews Genetics}, 16:146--158, 2015.

\bibitem{MP43}
W.~S. Mac~Culloch and W.~S. Pitts.
\newblock A logical calculus of the ideas immanent in nervous activity.
\newblock {\em Bull. Math Bio. Phys.}, 5:113--115, 1943.

\bibitem{PM95}
E.~Plathe, T.~Mestl, and S.W. Omholt.
\newblock Feedback loops, stability and multistationarity in dynamical systems.
\newblock {\em Journal of Biological Systems}, 3:569--577, 1995.

\bibitem{RRT08}
E.~Remy, P.~Ruet, and D.~Thieffry.
\newblock Graphic requirements for multistability and attractive cycles in a
  {B}oolean dynamical framework.
\newblock {\em Advances in Applied Mathematics}, 41(3):335--350, 2008.

\bibitem{R08}
A.~Richard.
\newblock An extension of a combinatorial fixed point theorem of {S}hih and
  {D}ong.
\newblock {\em Advances in Applied Mathematics}, 41(4):620--627, 2008.

\bibitem{R09}
A.~Richard.
\newblock Positive circuits and maximal number of fixed points in discrete
  dynamical systems.
\newblock {\em Discrete Applied Mathematics}, 157(15):3281--3288, 2009.

\bibitem{R10}
A.~Richard.
\newblock Negative circuits and sustained oscillations in asynchronous automata
  networks.
\newblock {\em Advances in Applied Mathematics}, 44(4):378--392, 2010.

\bibitem{R18}
A.~Richard.
\newblock Fixed points and connections between positive and negative cycles in
  boolean networks.
\newblock {\em Discrete Applied Mathematics}, 2018.

\bibitem{R19}
A.~Richard.
\newblock Positive and negative cycles in boolean networks.
\newblock {\em Journal of theoretical biology}, 463:67--76, 2019.

\bibitem{RC07}
A.~Richard and J.-P. Comet.
\newblock Necessary conditions for multistationarity in discrete dynamical
  systems.
\newblock {\em Discrete Applied Mathematics}, 155(18):2403--2413, 2007.

\bibitem{R86}
F.~Robert.
\newblock {\em Discrete iterations: a metric study}, volume~6 of {\em Series in
  Computational Mathematics}.
\newblock Springer, 1986.

\bibitem{R95}
F.~Robert.
\newblock {\em Les syst\`emes dynamiques discrets}, volume~19 of {\em
  Math\'ematiques et Applications}.
\newblock Springer, 1995.

\bibitem{S98}
E.H. Snoussi.
\newblock Necessary conditions for multistationarity and stable periodicity.
\newblock {\em Journal of Biological Systems}, 6:3--9, 1998.

\bibitem{S13}
S.~Soliman.
\newblock A stronger necessary condition for the multistationarity of chemical
  reaction networks.
\newblock {\em Bulletin of mathematical biology}, 75(11):2289--2303, 2013.

\bibitem{S03}
C.~Soul\'e.
\newblock Graphical requirements for multistationarity.
\newblock {\em Com{P}lex{U}s}, 1:123--133, 2003.

\bibitem{S06}
C.~Soul\'e.
\newblock Mathematical approaches to differentiation and gene regulation.
\newblock {\em C.R. Paris Biologies}, 329:13--20, 2006.

\bibitem{T73}
R.~Thomas.
\newblock {B}oolean formalization of genetic control circuits.
\newblock {\em Journal of Theoretical Biology}, 42(3):563--585, 1973.

\bibitem{T81}
R.~Thomas.
\newblock On the relation between the logical structure of systems and their
  ability to generate multiple steady states or sustained oscillations.
\newblock {\em Springer Series in Synergies 9}, pages 180--193, 1981.

\bibitem{TA90}
R.~Thomas and R.~d'Ari.
\newblock {\em Biological Feedback}.
\newblock CRC Press, 1990.

\bibitem{TK01a}
R.~Thomas and M.~Kaufman.
\newblock Multistationarity, the basis of cell differentiation and memory. {I}.
  structural conditions of multistationarity and other nontrivial behavior.
\newblock {\em Chaos: An Interdisciplinary Journal of Nonlinear Science},
  11(1):170--179, 2001.

\bibitem{TK01}
R.~Thomas and M.~Kaufman.
\newblock Multistationarity, the basis of cell differentiation and memory.
  {II}. {L}ogical analysis of regulatory networks in terms of feedback
  circuits.
\newblock {\em Chaos: An Interdisciplinary Journal of Nonlinear Science},
  11(1):180--195, 2001.

\end{thebibliography}

\end{document}